\documentclass[12pt]{amsart}
\usepackage{package}
\usepackage{subfiles}
\usepackage{hyperref}

\usepackage[
    backend=bibtex,
    giveninits=true,
    style=alphabetic,
    url=false
  ]{biblatex}
\title[Intertwining operators for $p$-adic covering groups]{Intertwining operators for representations of covering groups of reductive $p$-adic groups}
\date{\today}

\author{Janet Flikkema}
\address{\rm IMAPP, Radboud Universiteit, Heyendaalseweg 135, 6525AJ Nijmegen, the Netherlands}
\email{janet.flikkema@ru.nl}

\author{Maarten Solleveld}
\address{\rm IMAPP, Radboud Universiteit, Heyendaalseweg 135, 6525AJ Nijmegen, the Netherlands}
\email{m.solleveld@science.ru.nl}

\begin{document}

\begin{abstract}
    Let $G$ be a covering group of a reductive $p$-adic group. We study intertwining operators between parabolically induced representations of $G$ and prove that they satisfy certain adjointness relations. The Harish-Chandra $\mu$-function is defined as a composition of such intertwining operators for opposite parabolic subgroups of $G$. It can be seen as a complex rational function and we give an explicit formula for it in terms of poles and zeros. The adjointness of the intertwining operators is an important ingredient to prove the formula for the $\mu$-function. Moreover, we need to study the limit of the $\mu$-function at zero and infinity. To locate the poles of $\mu$, we construct a continuous family of Hermitian forms on a family of parabolically induced representations.
\end{abstract}

\maketitle

\tableofcontents

\section{Introduction}\label{sec:intro}

In this paper, we prove a formula for the Harish-Chandra $\mu$-function for covering groups of reductive $p$-adic groups. By a covering group $G$ of a reductive $p$-adic group $G'$, we mean a central extension of $G'$ by a finite abelian group, which is moreover a topological covering. An example of such a group is the metaplectic group, which is a double cover of the symplectic group. Many methods and results from the representation theory of reductive $p$-adic groups generalize to covering groups. For example, we have notions of parabolic induction and restriction functors, and there is a Bernstein decomposition, analogous to the reductive case. See \cite{FratilaPrasad} for an overview.

The Harish-Chandra $\mu$-function is defined using intertwining operators $J_{\bar{P}|{P}}$ and $J_{P|\bar{P}}$ between parabolically induced representations, for opposite parabolic subgroups $P$ and $\bar{P}$ of $G$. These intertwining operators are studied in \cite{Wal}; in fact, they can be considered more generally for parabolic subgroups $P$, $Q$ which have the same Levi subgroup, but are not necessarily opposite. Such intertwining operators also exist for covering groups, see  \cite{WenWeiLi}. An important ingredient in our work is an adjointness relation of these intertwining operators, which was stated but not proved in \cite{Wal}. Section \ref{sec:adjointness} of this paper gives a concrete proof of this adjointness property, using Bruhat-Tits theory. In section \ref{sec:polesandzeros}, we prove a formula for the Harish-Chandra $\mu$ function, as a complex rational function. The formula was already given for reductive groups in \cite{Sil}. Here, we provide an alternative proof of the result, which also works for covering groups.

To be able to more precisely state our results, let us give some more context. Let $M$ be a Levi subgroup of a covering group $G$ and let $P=MU_P$, $Q=MU_Q$ be parabolic subgroups of $G$ with Levi $M$. Let $\pi$ be any irreducible smooth $M$-representation, and let $I_P^G$ be the functor of normalized parabolic induction. Consider the intertwining operators 
\begin{align*}
        J_{Q|P}(\pi):& I_P^{{G}}(\pi) \rightarrow I_{{Q}}^{{G}}(\pi)\\
         &f \mapsto \left[g\mapsto \int_{({U_P}\cap {U_Q})\setminus {U_Q}}f(ug)\td u\right],
\end{align*}
 and $J_{P|Q}$ which is defined in a similar way. Let $\overline{\pi}^\vee$ denote the Hermitian dual of $\pi$. The $M$-invariant pairing $\langle \ , \ \rangle$ between $\pi$ and $\overline{\pi}^\vee$ gives rise to a $G$-invariant pairing $\langle \ , \ \rangle$ between $I_P^G(\pi)$ and $I_P^G (\overline{\pi}^\vee)$. Similarly, we get a pairing $\langle \ , \ \rangle$ between $I_{Q}^G(\pi)$ and $I_{Q}^G (\overline{\pi}^\vee)$. With respect to these pairings, we will prove the following result. 
 \begin{thm} The intertwining operators $J_{Q|P}(\pi)$ and $J_{P|Q}(\overline{\pi}^\vee)$ are adjoint, in the sense that for  $f_1\in I_P^G (\pi)$, $f_2\in I_{Q}^G(\overline{\pi}^\vee)$ we have
 \[
 \langle J_{Q|P}(\pi)f_1,f_2\rangle =\langle f_1, J_{P|Q}(\overline{\pi}^\vee)f_2\rangle.
 \]
 \end{thm}
There is a similar result for the contragredient representation $\pi^\vee$ instead of $\overline{\pi}^\vee$. In the special case where $\pi$ is Hermitian, we obtain the adjointness relation that was already stated in \cite{Wal} and which we will need to prove the result in Section \ref{sec:polesandzeros}.

Now, let us discuss the Harish-Chandra $\mu$-function. Consider an irreducible supercuspidal representation $\sigma$ of $M$. Let $P=MU$ be a parabolic subgroup of $G$ with Levi $M$, and let $\bar{P}=M\bar{U}$ be the parabolic subgroup of $G$ opposite to $P$. We have intertwining operators $J_{\bar{P}|P}(\sigma)$, $J_{P|\bar{P}}(\sigma)$ as above; they are integrals over $\bar{U}$ and $U$, respectively. 

The composition 
\[j(\sigma):=J_{P|\bar{P}}(\sigma) \circ J_{\bar{P}|P}(\sigma):I_P^{{G}}(\sigma)\rightarrow I_P^{{G}}(\sigma)
\] is a scalar, and does not depend on the choice of parabolic subgroup $P$ \cite{Wal}.
The Harish-Chandra $\mu$-function is, up to a positive real scalar, defined to be $j^{-1}$. This function $\mu$ plays a crucial role in the Plancherel formula for reductive $p$-adic groups \cite{Wal}, generalized to covering groups in \cite{WenWeiLi}. 

One may  view $\mu(\sigma\otimes\chi)$ as a complex rational function in the variable $\chi\in X_{\nr}(M)$, where $X_{\nr}(M)$ is the complex algebraic torus of unramified characters of $M$. The function $\mu$ then decomposes as a product $\mu=\prod_{\alpha}\mu_\alpha$, where $\alpha$ runs over the positive roots in the reduced root system corresponding to $M$ in $G$. Each $\mu_\alpha$ can be computed in a Levi subgroup $M_\alpha$ of $G$, which has $M$ as a maximal Levi subgroup. Moreover, $\mu_{\alpha}$ can be seen as a complex function in a single variable $z=\chi(h_\alpha^\vee)\in\mathbb{C}^\times$, where is $h_\alpha^\vee$ suitably chosen element of $M$.

\begin{thm} In the above context, there exists a unitary $\sigma_0\in\sigma\cdot X_{\nr}(M)$ such that $\mu_{\alpha}$ has the form
\begin{equation}\label{intro:formule} 
\mu_{\alpha}(\sigma_0\otimes \chi)=\mu_{\alpha}(\sigma_0,z) =c\cdot \frac{(1-z)(1-z^{-1})}{(1-qz)(1-qz^{-1})}\frac{(1+z)(1+z^{-1})}{(1+q'z)(1+q'z^{-1})},
\end{equation}
where $c\in\mathbb{R}_{>0}$ and $q,q'\in\mathbb{R}_{\geq 1}$.
\end{thm}
Both the zeros and the poles of $\mu_\alpha$ have clear representation-theoretic significance. If $\sigma\otimes \chi$ is a pole of $\mu_\alpha$, then $I_{P\cap M_\alpha}^{M_{\alpha}}(\sigma\otimes\chi)$ is reducible (but not decomposable). On the other hand, the set of roots $\alpha$ with $\mu_\alpha(\sigma\otimes \chi)=0$ forms a root system. It plays a role in the Knapp-Stein theory of intertwining operators \cite{WenWeiLi, silbergerknappstein}, which determines the decomposition of $I_P^G(\sigma\otimes \chi)$ in indecomposable representations. 

Let us give an outline of the proof of the theorem. By setting $\mu_{\alpha}(\sigma,0):=\lim_{z\rightarrow 0}\mu_{\alpha}(\sigma,z)$ and $\mu_{\alpha}(\sigma,\infty):=\lim_{z\rightarrow \infty}\mu_{\alpha}(\sigma,z)$, we can view $\mu_{\alpha}$ as a rational function on the projective curve $\mathbb{P}^1(\mathbb{C})$. Suppose $\mu_\alpha(\sigma,z)$ is not constant, then there exists a unitary $\sigma_0\in\sigma\cdot X_{\nr}(M)$ with $\mu_\alpha(\sigma_0)=0$. The only other possibility for a zero of $\mu_\alpha(\sigma_0,z)$ in $\mathbb{C}^\times$ is at $z=-1$. Every zero must be of order two, which corresponds to $J_{\bar{P}|P}$ and $J_{P|\bar{P}}$ each having a pole of order one at that point. So $\mu(\sigma_0,z)$ has at most two zeroes in $\mathbb{C}^\times$, located at $z=\pm 1$, each of order two. Moreover, we will show that $\mu(\sigma_0,z)$ is not zero at $z=0$ and $z=\infty$. We will also find two poles for every zero, using techniques concerning Hermitian and unitary representations. 

An application of our result is to give a more explicit description of Bernstein blocks in the category of smooth representations of $G$. It is known that Bernstein blocks can be realized as module categories of endomorphism algebras, see \cite{FratilaPrasad}. The work in \cite{solleveld2023endomorphismalgebrasheckealgebras} describes the structure of these endomorphism algebras for reductive $p$-adic groups. With the formula for the $\mu$-function for covering groups, and other results that are already known for covering groups, the results in \cite{solleveld2023endomorphismalgebrasheckealgebras} should generalize.\\

\textbf{Acknowledgement.}
The authors thank Caihua Luo for pointing out a problem in the first version of this paper.
We also thank Eric Opdam for his questions during the PhD defence of the first author, which
prompted us to expand the arguments in Paragraph \ref{par:zeros}.

\section{Notation}

$F$: non-archimedean local field\\
$G' = \mc G' (F)$: reductive group over $F$\\
$S' \subset G'$: maximal $F$-split torus\\
$M' \subset G'$: Levi subgroup containing $S'$\\
$P' = M' U_{P'}$: parabolic subgroup of $G'$\\
$\bar{P'} = M' U_{\bar{P'}}$: opposite parabolic subgroup\\
$p_G : G \to G'$: finite central covering\\
$M = p_G^{-1}(M')$: Levi subgroup of $G$\\
$P = p_G^{-1} (P'), \bar P = p_G^{-1} (\bar{P'})$: parabolic subgroups of $G$

Since $p_G$ splits over unipotent subgroups, $U_{P'}$ can be identified with
a subgroup $U_P \subset P$, and similarly $U_{\bar{P'}} \cong U_{\bar P} \subset \bar P$.
Then $P = M U_P$ and $\bar P = M U_{\bar P}$.\\
$K' = G'_x$: good maximal compact subgroup of $G'$, associated to a special vertex
$x$ in the apartment for $S'$ in the Bruhat--Tits building $\mc B (G')$ \\
$K = p_G^{-1} (K') = G_x$: good maximal compact subgroup of $G$\\
$\pi$: irreducible smooth $M$-representation on $\C$-vector space $V_\pi$\\
$I_P^G (\pi)$: normalized parabolic induction of $\pi$, on the vector space $I_P^G (V_\pi)$\\
$(\overline{\pi}^\vee, \overline{V_\pi}^\vee)$: Hermitian dual of $(\pi,V_\pi)$\\
$(\pi^\vee, V_{\pi}^\vee)$: contragredient of $(\pi, V_{\pi})$\\
$(\sigma,V_{\sigma})\in \Irr(M)_{\cusp}$: irreducible supercuspidal representation of $M$\\
$X_{\nr}(M)$: the torus of unramified characters of $M$\\
$\Irr(M)_{[M,\sigma]}=\sigma\cdot X_{\nr}(M)$: the inertial equivalence class of $\sigma$ for $M$

\section{Adjointness relations for intertwining operators}\label{sec:adjointness}

Let $G$ be a covering group of a reductive $p$-adic group and let $M$ be a Levi subgroup of $G$. In this section, we prove adjointness relations for the intertwining operators $J_{P|Q}$ and $J_{Q|P}$, where $P$ and $Q$ are parabolic subgroups of $G$ with Levi subgroup $M$. The adjointness relations will first be proved for opposite parabolic subgroups $P$ and $\bar{P}$, after which they will be generalized to arbitrary parabolic subgroups $P$ and $Q$ with the same Levi $M$.

Let $(\pi,V_{\pi})$ be an irreducible smooth representation of $M$, and denote by $(\overline{\pi}^\vee, \overline{V_{\pi}}^\vee)$ the Hermitian dual of $\pi$. The $M$-invariant sesquilinear pairing between $\pi$ and $\overline{\pi}^\vee$ gives rise
to a $G$-invariant sesquilinear pairing between $I_P^G (\pi)$ and $I_P^G (\overline{\pi}^\vee)$,
namely
\begin{equation}\label{eq:A.1}
\langle f, f' \rangle = \int_K \langle f(k), f' (k) \rangle \, \td k \qquad
f \in I_P^G (V_\pi), f' \in I_P^G (\overline{V_\pi}^\vee) .
\end{equation}
By the Iwasawa decomposition $G = P K$, $\int_K$ in \eqref{eq:A.1} is up to 
a positive scalar factor the same as $\int_{P \backslash G}$.
The pairing \eqref{eq:A.1} is nondegenerate, so provides an isomorphism 
$I_P^G (\overline{\pi}^\vee) \cong \overline{I_P^G (\pi)}^\vee$ .

There are intertwining operators 
\begin{equation}\label{eq:A.4}
J_{\bar P | P}(\pi) : I_P^G (\pi) \to I_{\bar P}^G (\pi) ,\quad J_{P | \bar P}
(\overline{\pi}^\vee) : I_{\bar P}^G (\overline{\pi}^\vee) \to I_P^G (\overline{\pi}^\vee) ,
\end{equation}
defined by
\begin{equation}\label{eq:intertwiningdef}
(J_{\bar P | P}(\pi) f) (g) = \int_{U_{\bar P}} f (\bar u g) \, \td \bar u ,\quad
(J_{P | \bar P}(\overline{\pi}^\vee) f') (g) = \int_{U_P} f' (u g) \, \td u .
\end{equation}
These arise via meromorphic continuation, as in \cite[\S IV.1]{Wal}.

We can build two $G$-invariant sesquilinear pairings
$I_P^G (V_\pi) \times I_{\bar P}^G (\overline{V_\pi}^\vee) \to \C$:
\begin{align}
\label{eq:A.2} & (f_1,f_2) \mapsto \langle J_{\bar P | P}(\pi) f_1, f_2 \rangle \\
\label{eq:A.3} & (f_1,f_2) \mapsto \langle f_1, J_{P | \bar P}(\overline{\pi}^\vee) f_2 \rangle, 
\end{align}
which are well-defined provided Condition \ref{cond:A.2}(i) below.
\begin{cond}\label{cond:A.2}
\begin{itemize}
\item[(i)] $J_{\bar P | P}(\pi)$ and $J_{P | \bar P}(\overline{\pi}^\vee)$ are regular, i.e. the
meromorphic continuation of (\ref{eq:intertwiningdef}) does not have a pole at $\pi$ or $\overline{\pi}^\vee$. 
\item[(ii)] $I_P^G (\pi)$ and $I_{\bar P}^G (\pi)$ are irreducible. This implies that 
$I_P^G (\overline{\pi}^\vee) \cong \overline{I_P^G (\pi)}^\vee$ and 
$I_{\bar P}^G (\overline{\pi}^\vee) \cong \overline{I_{\bar P}^G (\pi)}^\vee$ are also irreducible.
\end{itemize}
\end{cond}
The second condition holds for irreducible supercuspidal $M$-representations in general position by 
\cite[Th\'eor\`eme VI.8.5]{Ren}, and (more generally) for irreducible $M$-representations in general position 
by \cite[Th\'eor\`eme 3.2]{Sau}. Although these results are written for reductive $p$-adic groups,
the arguments apply also to finite covers like $G$. We elaborate on this below.

\begin{thm}\label{thm:Sau}
	Let $\pi$ be an irreducible representation of $M$. Let $\mathcal{V}$ be an open neighbourhood of $1$ in the variety $X_\nr(M)$ with the analytic topology. Then there exists a dense open subset $\mathcal{V}'$ of $\mathcal{V}$ such that $I_P^G(\pi\otimes\chi)$ is irreducible for all $\chi\in \mathcal{V}'$. 
\end{thm}
\begin{proof}
For reductive groups, this was proved in \cite[Th\'eor\`eme 3.2]{Sau}. We explain how the proof generalizes to covering groups. Most results of \cite{Sau} hold for locally compact unimodular groups, so in particular for finite central covering groups of reductive $p$-adic groups. The result in \cite[Th\'eor\`eme 3.2]{Sau} is specific for reductive $p$-adic groups and uses the following additional properties of representations of such groups: supercuspidal representations, parabolic induction, Jacquet restriction, Frobenius reciprocity, Bernstein's geometrical lemma \cite[Th\'eor\`eme IV.5.1]{Ren}, unramified characters and supercuspidal support. All of these concepts are worked out in \cite{Ren} and generalize to covering groups, as worked out in \cite{FratilaPrasad}. Therefore, \cite[Th\'eor\`eme 3.2]{Sau} generalizes to covering groups. 
\end{proof}

\begin{lem}\label{lem:A.1}
Assuming Condition \ref{cond:A.2}, the pairings \eqref{eq:A.2} and \eqref{eq:A.3} differ by a nonzero scalar factor.
\end{lem}
\begin{proof}
Any nonzero $G$-invariant sesquilinear pairing $V \times W \to \C$ defines a nonzero homomorphism 
of $G$-representations $\phi : V \to \overline{W}^\vee$, and conversely, by the formula
\[
(v,w) = \phi (v) (w) .
\]
If $V$ and $W$ are irreducible $G$-representations, then $\phi$ is unique up to scalars from
$\C^\times$, by Schur's lemma. The pairings $I_P^G (V_\pi) \times I_{\bar P}^G (\overline{V_\pi}^\vee) 
\to \C$ in \eqref{eq:A.2} and \eqref{eq:A.3} are nonzero because \eqref{eq:A.1} is nondegenerate and  
the operators \eqref{eq:A.4} are nonzero \cite[p. 283]{Wal}. Hence they differ only by a 
factor from $\C^\times$.
\end{proof}

We want to show that the pairings \eqref{eq:A.2} and \eqref{eq:A.3} are equal. This is claimed, but not
proven, for reductive $p$-adic groups in \cite[p. 287]{Wal}. It requires some preparation. 

Let $\Phi(G',S')$ be the root system of $G'$ and write it as 
\[
\Phi(G',S')=\Phi(M',S')\cup \Phi(U_P,S')\cup \Phi(U_{\bar P},S').
\]
Here, $\Phi(U_P,S')\subset \Phi(G',S')$ is the set of roots such that $U_P$ is generated by all root subgroups $U_\alpha$, where $\alpha\in \Phi(U_P,S')$. The subsets 
$\Phi(M', S')$ and $\Phi(U_{\bar P},S')$ are given similarly.

For each $r \in \R$, the valuated root datum of $G$ provides compact open subgroups $U_{\alpha,r}
\subset U_\alpha$. We normalize the valuations on the root subgroups $U_\alpha$ (of $G$ and of $G'$) so that 
\[
U_{\alpha,0} = U_\alpha \cap K = U_{\alpha,x} .
\]
By \cite[Lemma 7.3.11.(2)]{KaPr} there are well-defined subgroups
\[
U_{P,r} = \prod_{\alpha \in \Phi (U_P,S')} U_{\alpha,r} ,\quad
U_{\bar P,r} = \prod_{\alpha \in \Phi (U_{\bar P},S')} U_{\alpha,r} .
\]
For $r \geq 0$ we also have the subgroups $Z_G (S')_r = p_G^{-1} (Z_{G'}(S')_r)$ of $Z_G (S')$ 
and the Moy--Prasad groups \cite[\S 13.2]{KaPr}
\[
G'_{x,r} \subset G' \text{ and } G_{x,r} = p_G^{-1} (G'_{x,r}).
\]  
Recall that $p_G : G \to G'$ is a local homeomorphism and that the groups $G'_{x,r}$ with 
$r \in \R_{\geq 0}$ form a neighborhood basis of $e$ in $G'$ \cite[Proposition 13.2.5]{KaPr}. 
Restricting $p_G$ to a neighborhood of $e$ in $G$ on which it is a homeomorphism, we find
$r_x > 0$ and a homomorphic splitting $s_G : G'_{x,r_x} \to G$ of $p_G$. Then
\[
G_{x,r} = s_G (G'_{x,r}) \times \ker (p_G) \qquad \forall r \geq r_x.
\]
The advantage of this construction is that the subgroups $s_G (G'_{x,r})$ with $r \geq r_x$ form
a neighborhood basis of $e$ in $G$, in contrast with the $G_{x,r}$.
As $x$ lies in an apartment of $\mc B (M')$, we have $M' \cap K' = M'_x$ and $M \cap K = M_x$, 
and the groups $M'_{x,r}$ and $M_{x,r} = p_G^{-1}(M'_{x,r})$ are defined. 

Let $v_1 \in V_\pi, v_2 \in \overline{V_\pi}^\vee$, and pick $r \in \R_{>0}$ so that 
$s_G (M'_{x,r})$ fixes $v_1$ and $v_2$. We define $f_1 \in I_P^G (V_\pi)$ and 
$f_2 \in I_{\bar P}^G (\overline{V_\pi}^\vee)$ by the conditions
\begin{itemize}
\item $f_1 (e) = v_1$ and $f_2 (e) = v_2$,
\item supp$(f_1) = P U_{\bar P,r}$ and supp$(f_2) = \bar P U_{P,r}$,
\item $f_1$ is right $U_{\bar P,r}$-invariant and $f_2$ is right $U_{P,r}$-invariant.
\end{itemize}
We note that these conditions are not overdetermined because the multiplication maps
\begin{equation}\label{eq:A.9}
U_P \times M \times U_{\bar P} \to G \quad \text{and} \quad U_{\bar P} \times M \times U_P \to G
\end{equation}
are injective. That follows from the analogous statements for $G'$, which can be found for 
instance in \cite[Proposition 2.1.8.(3)]{CGP}.

We want to analyse
\begin{equation}\label{eq:A.5}
\langle J_{\bar P | P} (\pi) f_1, f_2 \rangle = \int_{U_{\bar P}} \int_K 
\langle f_1 (\bar u k), f_2 (k) \rangle \, \td k \td \bar u .
\end{equation}

\begin{lem}\label{lem:A.3}
The only nonzero contributions to \eqref{eq:A.5} come from \\ $k \in U_{\bar P,0} M_x U_{P,r}$
and $\bar u \in U_{\bar P,0}$ such that $\bar u k \in U_{\bar P,r} M_x U_{P,r}$.
\end{lem}
\begin{proof}
From the support of $f_2$ we see that we can only get nonzero contributions to \eqref{eq:A.5}
when $k \in K \cap \bar P U_{P,r}$. Write
\[
k = \bar u_1 m_1 u_2 \text{ with } \bar u_1 \in U_{\bar P}, m_1 \in M, u_2 \in U_{P,r}.
\]
As $U_{P,r} \subset K$, we have
\begin{equation}\label{eq:A.6}
\bar u_1 m_1 \in K \cap U_{\bar P} M = G_x \cap U_{\bar P} M = G_x \cap \bar P.
\end{equation}
By \cite[Proposition 8.3.1]{KaPr} there exists an integral model $\mc G_x'$ of $\mc G'$,
such that $\mc G_x' (\mf o_F) = G'_x$. Then $\bar{P'} \cap G'_x, U_{\bar{P'}} \cap G'_x$ and
$M' \cap G'_x$ determine $\mf o_F$-subgroup schemes $\bar{\mc P}_x', \mc U_{\bar{P'},x}'$
and $\mc M_x'$ of $\mc G_x'$. Since $\bar{\mc P} = \mc U_{\bar P} \rtimes \mc M'$ in
$\mc G'$ and $\mc G_x'$ is an integral model of $\mc G'$, also
\[
\bar{\mc P}_x' = \mc U_{\bar{P'},x}' \rtimes \mc M_x' .
\]
Taking $\mf o_F$-rational points, we obtain
\[
\bar{P'} \cap G'_x = \bar{\mc P}_x'  (\mf o_F) = \mc U_{\bar{P'},x}' (\mf o_F) \rtimes 
\mc M_x' (\mf o_F) = (U_{\bar P} \cap G'_x) \rtimes (M' \cap G'_x) = U_{\bar{P'},0} \rtimes M'_x .
\]
Applying $p_G^{-1}$, we find that
\[
\bar P \cap G_x = U_{\bar P,0} \rtimes M_x .
\]
Now \eqref{eq:A.6} says that $\bar u_1 \in U_{\bar P,0}$ and $m_1 \in M_x$, so
\begin{align*}
& k = \bar u_1 m_1 u_2 \in U_{\bar P,0} M_x U_{P,r}, \\
& \bar u k = \bar u \bar u_1 m_1 u_2 \in U_{\bar P} M_x U_{P,r} .
\end{align*}
The support of $f_1$ shows that $f_1 (\bar u k)$ can only be nonzero when 
$\bar u k \in U_P M U_{\bar P,r}$. We write 
\[
\bar u k = u_3 m_2 \bar u_5 \text{ with } u_3 \in U_P, m_2 \in M, \bar u_5 \in U_{\bar P,r}.
\]
Then
\begin{equation}\label{eq:A.7}
\bar u k u_2^{-1} =  \bar u \bar u_1 m_1 = u_3 m_2 \bar{u_5} u_2^{-1} \in U_{\bar P} M_x \cap
U_P M U_{\bar P,r} U_{P,r}. 
\end{equation}
The Iwahori decomposition of the Moy--Prasad group $G'_{x,r}$ \cite[\S 13.2]{KaPr} entails 
that the multiplication map
\begin{equation}\label{eq:A.8}
M_{x,r} \times U_{P,r} \times U_{\bar P,r} \to G_{x,r} \text{ is bijective.} 
\end{equation}
This also works with any other ordering of the three factors. Hence
\begin{multline*}
U_P M U_{\bar P,r} U_{P,r} = U_P M M_{x,r} U_{\bar P,r} U_{P,r} = \\
U_P M M_{x,r} U_{P,r} U_{\bar P,r} = U_P M U_{P,r} U_{\bar P,r} = U_P M U_{\bar P,r} .
\end{multline*}
This enables us to rewrite \eqref{eq:A.7} as 
\[
\bar u \bar u_1 m_1 = u_6 m_3 \bar u_7 \text{ with } u_6 \in U_P, m_3 \in M, \bar u_7 \in U_{\bar P,r} .
\]
Next we observe that
\[
u_6 m_3 = \bar u \bar u_1 m_1 \bar u_7^{-1} \in U_P M \cap U_{\bar P} M_x U_{\bar P,r} =
U_P M \cap M_x U_{\bar P} .
\]
By \eqref{eq:A.9}, the last intersection is just $M_x$. This forces $u_6 = 1, m_3 = m_1 \in M_x$
and 
\[
\bar u \bar u_1 = m_1 \bar u_7 m_1^{-1} \in m_1 U_{\bar P,r} m_1^{-1} = U_{\bar P,r} .
\]
As $\bar u_1 \in U_{\bar P,0} \supset U_{\bar P,r}$, that also forces $\bar u \in U_{\bar P,0}$.
\end{proof}

Assume that the Haar measure on $G_{x,r} = U_{\bar P,r} M_{x,r} U_{P,r}$ is the product of the Haar 
measures on the three factors. Now we are ready to evaluate \eqref{eq:A.5}. 

\begin{lem}\label{lem:A.4}
$\langle J_{\bar P | P}(\pi) f_1, f_2 \rangle$ equals $\mr{vol}(U_{\bar P,r}) 
\mr{vol}(U_{\bar P,0}) \mr{vol}(M_x) \mr{vol}(U_{P,r}) \langle v_1, v_2 \rangle$.
\end{lem}
\begin{proof}
By Lemma \ref{lem:A.3} we may replace the integral over $k \in K$ in \eqref{eq:A.5} by the
integral over $\bar u_1 m_1 u_2 \in U_{\bar P,0} \times M_x \times U_{P,r}$. Lemma \ref{lem:A.3}
also says that we may replace the integral over $\bar u \in U_{\bar P,0}$ by the integral over
$\bar u_2 = \bar u \bar u_1^{-1} \in U_{\bar P,r}$. Using the properties of $f_2$ we compute
\begin{align*}
\langle J_{\bar P | P}(\pi) f_1, f_2 \rangle & = \int_{U_{\bar P,r}} \int_{U_{\bar P,0}}
\int_{M_x} \int_{U_{P,r}} \langle f_1 (\bar u_2 m_1 u_2), f_2 (\bar u_1 m_1 u_2) \rangle
\, \td u_2 \td m_1 \td \bar u_1 \td \bar u_2 \\
& = \mr{vol}(U_{\bar P,0})\int_{U_{\bar P,r}} \int_{M_x} \int_{U_{P,r}} \langle f_1 (\bar u_2 m_1 u_2), 
f_2 (m_1) \rangle \, \td u_2 \td m_1 \td \bar u_2 .
\end{align*}
Since $M_x$ is compact and normalizes $U_{\bar P,r}$, we may exchange the order of $\bar u_2$ and
$m_1$. That gives
\begin{align*}
& \mr{vol}(U_{\bar P,0})\int_{U_{\bar P,r}} \int_{M_x} \int_{U_{P,r}} \langle f_1 (m_1 \bar u_2 u_2), 
f_2 (m_1) \rangle \, \td u_2 \td m_1 \td \bar u_2 = \\
& \mr{vol}(U_{\bar P,0})\int_{U_{\bar P,r}} \int_{M_x} \int_{U_{P,r}} \langle \pi (m_1) f_1 (\bar u_2 u_2), 
\overline{\pi}^\vee (m_1) f_2 (e) \rangle \, \td u_2 \td m_1 \td \bar u_2 = \\
& \mr{vol}(U_{\bar P,0}) \mr{vol}(M_x) \int_{U_{\bar P,r}} \int_{U_{P,r}} \langle f_1 (\bar u_2 u_2), 
f_2 (e) \rangle \, \td u_2 \td \bar u_2 .
\end{align*}
Here $\bar u_2 u_2 \in s_G (G'_{x,r})$, so by \eqref{eq:A.8} it can be written as 
\[
\bar u_2 u_2 = u_3 m_3 \bar u_3 \text{ with } u_3 \in U_{P,r}, m_3 \in s_G (M'_{x,r}), 
\bar u_3 \in U_{\bar P,r} .
\]
The construction of $f_1$ entails that
\[
f_1 (\bar u_2 u_2) = f_1 (u_3 m_3 \bar u_3) = f_1 (m_3) = \pi (m_3) f_1 (e) = \pi (m_3) v_1 = v_1 .
\]
Thus our integral simplifies to
\begin{multline*}
\mr{vol}(U_{\bar P,0}) \mr{vol}(M_x) \int_{U_{\bar P,r}} \int_{U_{P,r}} \langle v_1, 
f_2 (e) \rangle \, \td u_2 \td \bar u_2 = \\
\mr{vol}(U_{\bar P,0}) \mr{vol}(M_x) \mr{vol} (U_{\bar P,r}) \mr{vol} (U_{P,r}) \langle v_1, 
v_2 \rangle . \qedhere
\end{multline*}
\end{proof}

Completely analogous to Lemma \ref{lem:A.4}, one can show that
\begin{equation}\label{eq:A.10}
\langle f_1, J_{P | \bar P}(\overline{\pi}^\vee) f_2 \rangle = \mr{vol}(U_{P,0}) \mr{vol}(M_x) 
\mr{vol} (U_{\bar P,r}) \mr{vol} (U_{P,r}) \langle v_1, v_2 \rangle . 
\end{equation} 
For a better comparison, we need to normalize the Haar measures on $G$, $M$, $U_P$, $U_{\bar P}$. That
can be done by fixing the volume of one compact open subgroup in each of these groups:
\begin{align*}
\mr{vol}(M_x) &= 1, \mr{vol}(U_{P,0}) = 1, \mr{vol}(U_{\bar P,0}) = 1,\\
\mr{vol}(G_{x,r}) &= \mr{vol} (U_{\bar P,r}) \mr{vol}(M_{x,r}) \mr{vol}(U_{P,r}) .
\end{align*}
The last definition says that (locally) the Haar measure on $G$ is the product of the Haar
measures on $U_{\bar P}, M$ and $U_P$.

We will also formulate a version of the upcoming theorem for the pairing between
a representation $\pi$ and its contragredient $\pi^\vee$. By the same formula as in 
\eqref{eq:A.1}, the $M$-invariant bilinear pairing between $\pi$ and $\pi^\vee$ yields a 
$G$-invariant bilinear pairing between $I_P^G (\pi)$ and $I_P^G (\pi^\vee)$.

\begin{thm}\label{thm:A.5}
Let $\pi$ be any irreducible smooth $M$-representation.
\begin{enumerate}[(a)]
\item For any $f_1 \in I_P^G (V_\pi), f_2 \in I_{\bar P}^G (\overline{V_\pi}^\vee)$:
\[
\langle J_{\bar P | P}(\pi) f_1 , f_2 \rangle = 
\langle f_1, J_{P | \bar P}(\overline{\pi}^\vee) f_2 \rangle .
\]
\item For any $f_1 \in I_P^G (V_\pi), f \in I_{\bar P}^G (V_\pi^\vee)$:
\[
\langle J_{\bar P | P}(\pi) f_1 , f \rangle = \langle f_1, J_{P | \bar P}(\pi^\vee) f \rangle .
\]
\end{enumerate}
\end{thm}
\begin{rem} If $J_{\bar{P}|P}(\pi)$ is singular, then the two sides of the equalities in (a) and (b) of Theorem \ref{thm:A.5} may be infinite.
\end{rem}
\begin{proof}[Proof of Theorem \ref{thm:A.5}]
(a) First we assume that Condition \ref{cond:A.2} holds. By \\ Lemma \ref{lem:A.1}, it
suffices to find one pair $(f_1,f_2)$ such that $\langle J_{\bar P | P}(\pi) f_1 , f_2 \rangle$
equals \\ $\langle f_1, J_{P | \bar P}(\overline{\pi}^\vee) f_2 \rangle$ and is nonzero.
This is achieved by Lemma \ref{lem:A.4} and \eqref{eq:A.10} with $v_1,v_2$ such that
$\langle v_1, v_2 \rangle \neq 0$. That proves the theorem assuming Condition \ref{cond:A.2}.

Consider a family of irreducible smooth $M$-representations $\pi_0 \otimes \chi$, where
$\chi$ runs through the group $X_{\mr{nr}} (M)$ of unramified characters of $M$. By the rationality of $J_{\bar P | P}(\pi_0 \otimes \chi)$ and 
$J_{P | \bar P}(\overline{\pi \otimes \chi}^\vee)$ as functions of $\chi \in X_{\mr{nr}} (M)$ \cite[Th\'eor\`eme IV.1.1]{Wal} and by Theorem \ref{thm:Sau}, Condition \ref{cond:A.2} holds for $\pi_0 \otimes \chi$ 
with $\chi$ in some nonempty open subset of $X_{\mr{nr}} (M)$ with the analytic topology. In other words, we have
already proven the theorem for all $\pi$ in a nonempty open subset of the variety
\[
\{ \pi_0 \otimes \chi : \chi \in X_{\mr{nr}} (M) \}
\]
with the analytic topology. Consider a $\chi_0$ for which Condition \ref{cond:A.2} fails, and let
$f_3 \in I_P^G (\pi \otimes \chi_0)$ and $f_4 \in I_{\bar P}^G (\overline{\pi_0 \otimes \chi_0}^\vee)$.
Via the linear bijections
\[
I_P^G (V_{\pi_0 \otimes \chi_0}) \cong I_{P \cap K}^K (V_{\pi_0 \otimes \chi_0}) \cong 
I_{P \cap K}^K (V_{\pi_0 \otimes \chi}) \cong I_P^G (V_{\pi_0 \otimes \chi})
\]
we can regard $f_3$ also as an element of $I_P^G (V_{\pi_0 \otimes \chi})$ for any
$\chi \in X_{\mr{nr}}(M)$. Similarly $f_4 \in I_{\bar P}^G (\overline{\pi_0 \otimes \chi}^\vee)$.
Now the expressions 
\[
\langle J_{\bar P | P}(\pi_0 \otimes \chi) f_3, f_4 \rangle \quad \text{and} \quad
\langle f_3, J_{P | \bar P}(\overline{\pi_0 \otimes \chi}^\vee) f_4 \rangle 
\]
are rational functions of $\chi \in X_{\mr{nr}}(M)$. They coincide on an analytic open, hence Zariski-dense subset of
$X_{\mr{nr}}(M)$, so they are equal.\\
(b) This can be shown in the same way as part (a), using $\pi^\vee$ instead of $\overline{\pi}^\vee$
in all the previous arguments. 
\end{proof}

Theorem \ref{thm:A.5}.a says that $J_{\bar P | P} (\pi) : I_P^G (\pi) \to I_{\bar P}^G (\pi)$
is the adjoint operator of 
\[
J_{P | \bar P}(\overline{\pi}^\vee) : \overline{I_{\bar P}^G (\pi)}^\vee \cong 
I_{\bar P}^G (\overline{\pi}^\vee) \to I_P^G (\overline{\pi}^\vee) \cong \overline{I_P^G (\pi)}^\vee .
\]
In the special case that $\pi$ is Hermitian (e.g. unitary), we may replace $\overline{\pi}^\vee$
by the isomorphic representation $\pi$. Then Theorem \ref{thm:A.5}.a says that
\begin{equation}\label{eq:A.11}
J_{\bar P | P} (\pi) : I_P^G (\pi) \to I_{\bar P}^G (\pi) 
\text{ is the adjoint of } J_{P | \bar P}(\pi) : I_{\bar P}^G (\pi) \to I_P^G (\pi) ,
\end{equation}
with respect to the $G$-invariant Hermitian forms \eqref{eq:A.1} on $I_P^G (V_\pi)$ and on
$I_{\bar P}^G (V_\pi)$.

Theorem \ref{thm:A.5}.b is mentioned (for reductive $p$-adic groups) in \cite[p. 283]{Wal}, 
but unfortunately without a proof. We note that in the versions of Theorem \ref{thm:A.5} in
\cite{Wal}, the two parabolic subgroups share the same Levi factor, but they need not be
opposite. We want to extend Theorem \ref{thm:A.5} to that generality, for our covering group $G$.

Let $Z(M')_s$ be the maximal $F$-split torus in the centre of $M'$. We say that a root of $Z(M')_s$
is a character of $Z(M')_s$ appearing in the adjoint action on $\mr{Lie}(G) = \mr{Lie}(G')$, and
that a root $\alpha$ is reduced if $q\alpha$ is not a root for $q \in (0,1)$.
Let $Q = M U_Q$ be another parabolic subgroup of $G$ with Levi factor $M$. 
As in \cite[p. 279]{Wal}, we define the distance between $P$ and $Q$ as 
\[
d (P,Q) = \text{ number of reduced roots of } Z(M')_s \text{ appearing in Lie}(U_{\bar P} \cap U_Q) .
\]
In this more general context the definition of $J_{Q | P}(\pi)$ becomes
\[
(J_{Q | P}(\pi) f) (g) = \int_{(U_P \cap U_Q) \backslash U_Q} f (u g) \, \td u .
\]

\begin{prop}\label{prop:A.6}
Theorem \ref{thm:A.5} holds also with $Q$ instead of $\bar P$, that is, it holds for
any two parabolic subgroups with the same Levi factor.
\end{prop}
\begin{proof}
We can find a sequence of parabolic subgroups 
\[
P = P_0, P_1, \ldots, P_d = Q
\]
such that $d (P_i,P_j) = i-j$ for all $i \geq j$. According to \cite[p. 283]{Wal}, 
\[
J_{Q | P}(\pi) = J_{P_d | P_{d-1}}(\pi) \circ \cdots \circ J_{P_2 | P_1}(\pi) \circ J_{P_1 | P_0}(\pi). 
\]
Hence it suffices to prove the proposition for two parabolic subgroups of distance 1. In other
words, we only need to consider the cases with $d(P,Q) = 1$. 

Now $Z(M')_s$ acts on Lie$(U_{\bar P} \cap U_Q)$ by multiples of a unique reduced root $\alpha$.
We put $T'_\alpha = \ker (\alpha)^\circ$, a codimension one subtorus of $Z(M')_s$. Consider the
Levi subgroup $M'_\alpha := Z_{G'}(T'_\alpha)$ of $G'$. It is generated by 
\[
M' \cup (U_{\bar P} \cap U_Q) \cup (U_P \cap U_{\bar Q})
\]
and contains $M'$ as a maximal Levi subgroup. Then $M_\alpha := p_G^{-1} (M'_\alpha)$ is a Levi
subgroup of $G$ with $M$ as maximal Levi subgroup. Since 
\[
(U_P \cap U_Q) \backslash U_Q \cong U_{\bar P} \cap U_Q 
\]
and $U_{\bar P} \cap U_Q \subset M_\alpha$, both $J_{Q | P}(\pi)$ and 
$J_{Q \cap M_\alpha | P \cap M_\alpha}(\pi)$ are given by integration over $U_{\bar P} \cap U_Q$. 
It follows that $J_{Q | P}(\pi)$ can be obtained from $J_{Q \cap M_\alpha | P \cap M_\alpha}(\pi)$ 
by applying the functor $I_{P M_\alpha}^G = I_{Q M_\alpha}^G$.

In $M_\alpha$ there are only two parabolic subgroups with Levi factor $M$, namely $P \cap M_\alpha$
and $Q \cap M_\alpha$, and they are opposite. Thus Theorem \ref{thm:A.5} applies to
$J_{Q \cap M_\alpha | P \cap M_\alpha}(\pi)$, for any irreducible smooth $M$-representation $\pi$.
We note that the group $K \cap M_\alpha = M_{\alpha,x}$ is the good maximal compact subgroup of 
$M_\alpha$ associated to the special vertex $x$.
For $f_1 \in I_P^G (V_\pi)$ and $f_2 \in I_Q^G (\overline{V_\pi}^\vee)$ we compute
\begin{align}
\nonumber \langle J_{Q | P} (\pi) f_1, f_2 \rangle & 
= \int_K \langle (J_{Q | P} (\pi) f_1)(k), f_2 (k) \rangle \, \td k \\
\label{eq:A.12} & = \int_{K \cap M_\alpha \backslash K} \int_{K \cap M_\alpha} \langle 
(J_{Q | P} (\pi) f_1)(k_\alpha k), f_2 (k_\alpha k) \rangle \, \td k_\alpha \td k \\
\nonumber & = \int_{K \cap M_\alpha \backslash K} \int_{K \cap M_\alpha} \langle (J_{Q | P} (\pi) 
\pi (k) f_1) (k_\alpha), (\overline{\pi}^\vee (k) f_2) (k_\alpha) \rangle \, \td k_\alpha \td k .
\end{align}
The restriction of $\pi (k) f_1 : G \to V_\pi$ to $M_\alpha$ can be regarded as an element of
$I_{P \cap M_\alpha}^{M_\alpha}(V_\pi)$, and similarly for $\overline{\pi}^\vee (k) f_2$.
Then Theorem \ref{thm:A.5}.a says that \eqref{eq:A.12} equals 
\begin{equation}\label{eq:A.13}
\begin{aligned}
& \int_{K \cap M_\alpha \backslash K} \langle J_{Q \cap M_\alpha | P \cap M_\alpha}(\pi)
(\pi (k) f_1) |_{M_\alpha}, (\overline{\pi}^\vee (k) f_2) |_{M_\alpha} \rangle \, \td k \; = \\
& \int_{K \cap M_\alpha \backslash K} \langle (\pi (k) f_1) |_{M_\alpha}, J_{P \cap M_\alpha | 
Q \cap M_\alpha}(\overline{\pi}^\vee) (\overline{\pi}^\vee (k) f_2) |_{M_\alpha} \rangle \, \td k \; = \\
& \int_{K \cap M_\alpha \backslash K} \int_{K \cap M_\alpha} \langle f_1 (k_\alpha k), (J_{P \cap
M_\alpha | Q \cap M_\alpha}(\overline{\pi}^\vee) f_2) (k_\alpha k) \rangle \, \td k_\alpha \td k \; = \\
& \int_K \langle f_1 (k), (J_{P | Q} (\overline{\pi}^\vee) f_2)(k) \rangle \, \td k \; = \; 
\langle f_1, J_{P | Q} (\overline{\pi}^\vee) f_2 \rangle .
\end{aligned}
\end{equation}
The combination of \eqref{eq:A.12} and \eqref{eq:A.13} shows that Theorem \ref{thm:A.5}.a holds
for $J_{Q | P}(\pi)$ and $J_{P | Q}(\overline{\pi}^\vee)$. 
The same argument applies with $\pi^\vee$ instead of $\overline{\pi}^\vee$, and that generalizes
Theorem \ref{thm:A.5}.b.
\end{proof}

\section{Poles and zeros of the Harish-Chandra \texorpdfstring{$\mu$}{mu}-function}\label{sec:polesandzeros}

In this section, we discuss the Harish-Chandra $\mu$-function. By \cite[Lemme V.2.1]{Wal}, it decomposes as a product $\mu=\prod_{\alpha}\mu_\alpha$, where $\alpha$ runs over the positive roots in the reduced root system corresponding to $M$ in $G$. Each $\mu_\alpha$ can be computed in a Levi subgroup $M_\alpha$ of $G$, which has $M$ as a maximal Levi subgroup.

Thus from now on, we assume that $M$ is a maximal Levi subgroup of $G$, which leads to $[N_G(M):M]\leq 2$. Let $P=MU_{P}$ be a parabolic subgroup with Levi $M$ and let $\bar{P}=MU_{\bar P}$ be its opposite parabolic subgroup. We consider the intertwining operators  
\begin{equation*}
J_{\bar P | P}(\sigma\otimes\chi) : I_P^G (\sigma\otimes\chi) \to I_{\bar P}^G (\sigma\otimes\chi) ,\quad J_{P | \bar P}
(\sigma\otimes \chi) : I_{\bar P}^G (\sigma\otimes \chi) \to I_P^G (\sigma \otimes\chi) ,
\end{equation*}
for $\sigma\in \Irr_\cusp (M)$ and $\chi\in X_\nr(M)$. They are defined by
\[
(J_{\bar P | P}(\sigma\otimes\chi) f) (g) = \int_{U_{\bar P}} f (\bar u g) \, \td \bar u ,\quad
(J_{P | \bar P}(\sigma\otimes\chi) f') (g) = \int_{U_P} f' (u g) \, \td  u .
\]
Their composition
\[
j(\sigma\otimes \chi):=J_{P | \bar P}(\sigma\otimes\chi)\circ J_{\bar P| P}(\sigma\otimes \chi)
\]
is a rational function in the variable $\chi\in X_{\nr}(M)$, see \cite[\S IV.3]{Wal}. Up to a positive real scalar, the Harish-Chandra $\mu$-function is $j(\sigma\otimes\chi)^{-1}$ \cite[\S V.2]{Wal}, and we want to analyze its poles and zeros.

\subsection{\texorpdfstring{$\mu$}{mu} as a rational function}

Let us explain how to view $\mu(\sigma\otimes \chi)$ as a rational function on $\mathbb{C}^\times$. We refer to \cite{endomorphismalgebras} (and the correction in \cite{solleveld2023endomorphismalgebrasheckealgebras}).
Let
\[
X_\nr(M,\sigma):=\set{\chi\in X_\nr(M):\sigma\otimes \chi\cong \sigma}.
\]
Note that $X_\nr(M,\tau)=X_\nr(M,\sigma)$ for any $\tau\in \Irr(M)_{[M,\sigma]}$. There are subgroups of $M$,
\[
M^1=\cap_{\chi\in X_\nr(M)} \ker(\chi), \ M^2_\sigma=\cap_{\chi\in X_\nr(M,\sigma)}\ker (\chi).
\]
The subgroup $M\subset G$ is the inverse image of a maximal Levi subgroup $M'$ of $G'$. Let $Z(M')_s$ be the maximal split torus in the center of $M'$. Then there is only one positive reduced root $\alpha$ coming from the adjoint representation of $Z(M')_s$ on the Lie algebra of $G'$. Thus we have $(M_{\sigma}^2\cap G^1)/M^1\cong\mathbb{Z}$. By the following result, we can view $\mu$ as a rational function on $\mathbb{C}^\times$. 

\begin{lem}\label{lem:mudependsongenerator}
	In the context above, let $m$ be a generator of $(M^2_\sigma\cap G^1)/M^1\cong \mathbb{Z}$. Then $\mu(\sigma\otimes\chi)$ only depends on $\chi(m)$.
\end{lem}
\begin{proof}
	The proof is inspired by considerations in \cite[\S 2]{parameters} and we reduce to the derived subgroup. Let $\mathcal{G}'_\der$ be the derived subgroup of $\mathcal{G}'$, $G'_\der:=\mathcal{G}'_\der(F)$ its group of $F$-points, and let $G_\der$ be the inverse image of $G'_\der$ in $G$. The inclusion map $\mathcal{G}'_\der\rightarrow \mathcal{G}'$ is a central isogeny. It induces an inclusion $G'_\der\rightarrow G'$ between the groups of $F$-points, which lifts to an inclusion $G_\der \rightarrow G$ between covering groups. The Levi subgroups of $G$ are in one-to-one correspondence with the Levi subgroups of $G_\der$, via the bijection $M\mapsto \tilde M:=M\cap G_\der$.
	
	We have a restriction functor $-|_{\tilde{M}}:\Rep(M)\rightarrow \Rep(\tilde M)$, restricting to an algebraic group homomorphism $X_\nr(M)\rightarrow X_\nr(\tilde M)$. For $\sigma$ an irreducible supercuspidal representation of $M$, the restriction $\sigma|_{\tilde{M}}$ is a finite direct sum of irreducible supercuspidal $\tilde M$-representations. This was proved for reductive groups in \cite{isogenyrestrictions}, but the arguments generalize to our situation of covering groups.
	
	Since $M$ is a maximal Levi subgroup of $G$, $\tilde M$ is a maximal Levi subgroup of $G_\der$. Let $\sigma\in \Irr_\cusp(M)$ and let $\tilde{\sigma}\in\Irr_\cusp(\tilde M)$ be an irreducible constituent of $\sigma|_{\tilde M}$. We claim that there exists $\tilde{c}\in\mathbb{C}^\times$ such that
	\begin{equation}\label{eq:muderivedgp}
	\mu^{M}(\sigma\otimes\chi)=\tilde c\cdot \mu^{\tilde M}(\tilde{\sigma}\otimes \chi|_{\tilde M}).
	\end{equation}
	For reductive groups, this was proved in \cite[Proposition 2.2]{parameters}, however the proof refers to the explicit formula for the $\mu$-function which is known for reductive groups, but at this point we do not yet have it available for covering groups. Therefore, we have to be a bit careful. Since $\mu^M$ and $\mu^{\tilde M}$ are both rational functions, it suffices to show that they have the same poles and the same zeros. The poles have been taken care of in the proof of \cite[Proposition 2.2]{parameters}; the proof generalizes to covering groups, since it only uses the definition $\mu(\sigma\otimes\chi):=j(\sigma\otimes\chi)^{-1}$ (up to a scalar), rather than the explicit formula for $\mu$. Showing that the $\mu$-functions have the same zeros can be done in a similar way. So the claim in Equation (\ref{eq:muderivedgp}) holds.
	
	Thus $\mu(\sigma\otimes\chi)$ only depends on $\chi|_{\tilde M}$ (recall that $\tilde M=M\cap G_\der$). Since $\mu(\sigma\otimes\chi)$ is invariant under $X_\nr(M,\sigma)$, it follows that $\mu(\sigma\otimes \chi)$ only depends on $\chi|_{M^2_\sigma\cap G_\der}$. There is a canonical homomorphism between the quotients
	\[
	\mathbb{Z}\cong(M^2_\sigma\cap G_\der)/ (M^1\cap G_\der)\rightarrow (M^2_\sigma\cap G^1)/M^1\cong \mathbb{Z}.
	\]
	By the considerations above, $\mu(\sigma\otimes\chi)$ only depends on $\chi(m')$, for $m'$ a generator of $(M^2_\sigma\cap G_\der)/ (M^1\cap G_\der)$. The image of $m'$ in $(M^2_\sigma\cap G_1)/M^1$ is an integer power of $m$, where $m$ is a generator of $(M^2_\sigma\cap G_1)/M^1$. So $\chi(m')$ is an integer power of $\chi(m)$. Since $\mu(\sigma\otimes\chi)$ only depends on $\chi(m')$, it only depends on $\chi(m)$. 
\end{proof}

Recall that since $M$ is maximal, we consider only one root $\alpha$. Let $h_\alpha^\vee$ be the unique generator of $(M_{\sigma}^2\cap G^1)/M^1\cong\mathbb{Z}$ such that $|\alpha(h_\alpha^\vee)|_F>1$. Then by Lemma \ref{lem:mudependsongenerator}, we may view $\mu(\sigma\otimes\chi)=\mu(\sigma,z)$ as a function in the variable $z:=\chi(h_\alpha^\vee)\in \mathbb{C}^\times$. Suppose $\tau=\sigma\otimes \chi_{\tau}\in \Irr(M)_{[M,\sigma]}$ is a different choice of representative in the inertial equivalence class of $\sigma$ for $M$. Then for $\chi\in X_{\nr}(M)$, $\mu(\tau\otimes\chi)=\mu(\sigma\otimes\chi_{\tau}\chi)$. Hence
\begin{equation}\label{eq:mutau}
    \mu(\tau,z)=\mu(\sigma,q_{\tau}z), \textnormal{ where } q_\tau=\chi_{\tau}(h_\alpha^\vee)\in \mathbb{C}^\times.
\end{equation}

For $G'$, it was shown in \cite[Theorem 1.6]{Sil} that there exists a unitary $\sigma_0\in \Irr(M')_{[M',\sigma]}$ such that as a function in $z=\chi(h_\alpha^\vee)$, $\mu=\mu_\alpha$ has the form 
\begin{equation}\label{eqn:mu}
\mu(\sigma_0\otimes \chi)=\mu(\sigma_0,z) =c_{\mu}\cdot \frac{(1-z)(1-z^{-1})}{(1-qz)(1-qz^{-1})}\frac{(1+z)(1+z^{-1})}{(1+q'z)(1+q'z^{-1})},
\end{equation}
where $c_{\mu}\in\mathbb{R}_{>0}$ and $q,q'\in \mathbb{R}_{\geq 1}$. The possibility that $\mu$ is independent of $\chi$ is included in this formula; this happens when $q=q'=1$. 

Our goal is to give an alternative proof of this result, which generalizes to covering groups. The idea is as follows. Since $\mu(\sigma\otimes\chi)=\mu(\sigma,z)$ is a rational function on $\mathbb{C}^\times$, we may also view it as a rational function on the projective line $\mathbb{P}^1(\mathbb{C})$, by considering the limits $\mu(\sigma,\infty):=\lim_{z\rightarrow\infty}\mu(\sigma,z)$ and $\mu(\sigma,0):=\lim_{z\rightarrow 0}\mu(\sigma,z)$. If $\mu(\sigma,z)$ is not constant, then there exists a unitary $\sigma_0\in \Irr(M)_{[M,\sigma]}$ with $\mu(\sigma_0)=0$. The only other possibility for a zero of $\mu(\sigma_0,z)$ in $\mathbb{C}^\times$ will be at $z=-1$ and all zeros are of multiplicity $2$. Using techniques concerning Hermitian and unitary representations (see \cite[\S IV]{Ren} for definitions), we will then locate two poles in $\mathbb{R}^\times$ for every zero of $\mu(\sigma_0,z)$. Moreover, we will show that $\mu(\sigma_0,z)$ can not be zero at $z=\infty$ and $z=0$. Combining all these results, we obtain the formula as in (\ref{eqn:mu}). 

The result below is similar to the reductive case and will help us to give a description of the supercuspidal representations of $M$.

\begin{lem}\label{lem:subgpfiniteindex}
    Let $p_G:G\rightarrow G'$ be a finite central covering. Let $G^1=p_G^{-1}(G'^1)$, where $G'^1$ is the subgroup of $G'$ generated by all compact subgroups. Then $G^1\cdot Z(G)$ has finite index in $G$.
\end{lem}
\begin{proof}
    Let $g\in p_G^{-1}(Z(G'))$. Consider the commutator map \[[g, \cdot ] :G\rightarrow\ker(p_G)\subset Z(G), \ x\mapsto gxg^{-1}x^{-1}.\] Let $k=|\ker (p_G)|$. Then since $gxg^{-1}=[g,x]x$, we obtain
    \[
    g^kxg^{-k}x^{-1}=g^{k-1}[g,x]xg^{1-k}x^{-1}=\ldots=[g,x]^kxx^{-1}=[g,x]^k=1.
    \]
    Therefore, $g^k\in Z(G)$ for all $g\in p_G^{-1}(Z(G'))$. Recall that $G'/G'^1$ is a finitely generated free abelian group and that $Z(G')/(Z(G')\cap G'^1)$ is a subgroup of finite index in $G'/G'^1$. Hence the group of $k$-th powers $Z(G')^{(k)}/(Z(G')^{(k)}\cap G'^1)$ also has finite index in $G'/G'^1$. Applying $p_G^{-1}$, it follows that $G^1\cdot Z(G)$ has finite index in $G$.
\end{proof}

Let $\Irr(M)^{\unit}_{\cusp}$ be the set of unitary supercuspidal representations of $M$.

\begin{lem}\label{lem:supercuspidalunitary}
    There is a bijection
    \[
    \Irr(M)_\cusp^{\unit}\times \Hom(M,\mathbb{R}_{>0})\rightarrow \Irr(M)_\cusp.
    \]
\end{lem}
\begin{proof} 
First, let us show that there is a bijection 
\[
\Hom(Z(M),\mathbb{R}_{>0})\xrightarrow{\sim} \Hom(M,\mathbb{R}_{>0}).
\]
For this, we need to show that every $\chi\in\Hom(Z(M),\mathbb{R}_{>0})$ extends uniquely to a character in $\Hom(M,\mathbb{R}_{>0})$. Since $Z(M)$ is abelian, it has a maximal compact subgroup $Z(M)_{\cpt}$. Note that $\chi$ factors through $Z(M)/Z(M)_{\cpt}$, because the image of $Z(M)_{\cpt}$ must be a compact subgroup of $\mathbb{R}_{>0}$, i.e. $\set{1}$. By Lemma \ref{lem:subgpfiniteindex}, $M^1Z(M)$ has finite index in $M$. Hence, \[
Z(M)/Z(M)_{\cpt}\hookrightarrow M/M^1\cong \mathbb{Z}^r,
\]
is a subgroup of finite index. So $Z(M)/Z(M)_{\cpt}\cong \mathbb{Z}^r$ is a free abelian subgroup of the same rank. So for $x\in M/M^1$, let $n$ be a positive integer such that $n\cdot x \in Z(M)/Z(M)_{\cpt}$. Define $\chi(x):=\chi(n\cdot x)^{1/n}$, which is well-defined because $\mathbb{R}_{>0}$ has unique $n$-th roots and it does not depend on the choice of $n$.  It is a group homomorphism, since 
\begin{align*}
\chi(n\cdot x)^{\frac{1}{n}}\chi(m\cdot y)^{\frac{1}{m}}&=\chi(nm\cdot x)^{\frac{1}{nm}}\chi(nm\cdot y)^{\frac{1}{nm}}=(\chi(nm\cdot x)\chi(nm\cdot y))^{\frac{1}{nm}}\\
&=\chi((nm\cdot x)(nm\cdot y))^{\frac{1}{nm}}=\chi(nm\cdot xy)^{\frac{1}{nm}}.
\end{align*}
This describes how $\chi:Z(M)/Z(M)_{\cpt}\rightarrow \mathbb{R}_{>0}$ extends uniquely to a character $\chi:M/M^1\rightarrow \mathbb{R}_{>0}$.
Now we can explicitly describe the bijection. It is given by
\[
    \Irr(M)_\cusp^{\unit}\times \Hom(M,\mathbb{R}_{>0})\rightarrow \Irr(M)_\cusp, (\pi,\chi)\mapsto \pi\otimes \chi,
    \]
and its inverse is given by $\sigma\mapsto (\sigma\otimes |\cc(\sigma)|^{-1}, |\cc(\sigma)|)$, where $\cc(\sigma)$ is the central character of $\sigma$, and we extend $|cc(\sigma)|$ to a character of $M$ as above. By \cite[116]{Ren}, \cite[Theorem 6.2]{FratilaPrasad}, a smooth irreducible representation $\pi$ of $M$ is supercuspidal and unitary if and only if it is supercuspidal and its central character $\cc(\pi)$ is unitary. That explains why the representation $\sigma\otimes |cc(\sigma)|^{-1}$ is unitary.
\end{proof}

\begin{rem}
Although Lemma \ref{lem:supercuspidalunitary} is written for supercuspidal representations, it holds more generally for representations that are essentially square integrable, see \cite[\S IV.3.2]{Ren} for the definition.
\end{rem}

By Lemma \ref{lem:supercuspidalunitary}, there exists a unitary $\tilde\sigma\in \Irr(M)_{[M,\sigma]}$. Knowing $\mu(\tilde\sigma,z)$ is equivalent to knowing $\mu(\sigma,z)$, by (\ref{eq:mutau}).

Suppose that $\sigma$ is unitary, then by \cite[\S IV.2.3]{Ren}, the parabolically induced representation $I_P^G(\sigma)$ is unitary as well. We can use the adjointness relation from Section \ref{sec:adjointness} to show that $\mu$ must be real and non-negative on the unit circle.

\begin{lem}\label{lem:unitcircle}
Let $\sigma$ be unitary. Then the function $\mu(\sigma, z)$ has all its values in $\mathbb{R}_{\geq 0}$ for $z$ on the unit circle $S^1=\set{e^{i\varphi}\ |\ \varphi\in \R}\subset\C^\times$.
\end{lem}
\begin{proof}
    We claim that for $z=e^{i\varphi}\in S^1$, there always exists a unitary character $\chi_z\in X_\nr(M)$ such that $\chi_z(h_\alpha^\vee)=z$. Indeed, recall that $X_{\nr}(M)=\Hom(M/M^1,\mathbb{C}^\times)$. Moreover, $M/M^1$ is a free abelian group of finite rank, and $h_\alpha^\vee$ is a generator of $M_\sigma^2\cap G^1/M^1$, which is a free rank $1$ subgroup of $M/M^1$. Therefore, we can choose a basis of $M/M^1$, such that $h_\alpha^\vee$ is an integer multiple of one of these basis elements. So write $h_\alpha^\vee=n\cdot b$ with $n\in \mathbb{Z}$ and $b$ a basis element of $M/M^1$. Then we can define $\chi_z(b)=e^{i\varphi/n}$, and let $\chi_z$ be trivial on the other basis elements of $M/M^1$. This gives us a unitary character with 
    \[
    \chi_z(h_\alpha^\vee)=\chi_z(n\cdot b)=\chi_z(b)^n=(z^{1/n})^n =z.
    \]
    Since we chose $\sigma$ to be unitary, $\sigma\otimes\chi_z$ is unitary as well, and the adjunction formula from Theorem \ref{thm:A.5}.a applies with $\overline{\sigma\otimes\chi_z}^\vee$ replaced by $\sigma\otimes\chi_z$. Since $\sigma\otimes \chi_z$ is unitary, we may assume that the induced pairing on $I_P^G(\sigma\otimes\chi_z)$ (see (\ref{eq:A.1})) is positive definite.
    So for $f\in I_P^G(\sigma\otimes \chi_z)$, we have 
    \begin{align}\label{eq:unitcircle}
    j(\sigma\otimes\chi_z)\langle f, f\rangle &=\langle f, J_{P | \bar P}(\sigma\otimes\chi_z)\circ J_{\bar P| P}(\sigma\otimes \chi_z) f\rangle \\
    &\nonumber = \langle J_{\bar P| P}(\sigma\otimes \chi_z) f, J_{\bar P| P}(\sigma\otimes \chi_z) f\rangle\in\mathbb{R}_{\geq 0}\cup\{\infty\}.
    \end{align}
 Since $J_{P|\bar P}(\sigma\otimes\chi_z)$ is nonzero, (\ref{eq:unitcircle}) shows that $j(\sigma\otimes\chi_z)\not=0$. Thus $z\mapsto j(\sigma\otimes\chi_z)$ is a continuous map from $S^1$ to $\mathbb{R}_{>0}\cup \set{\infty}$. We obtain that $\mu(\sigma,z)=\mu(\sigma\otimes \chi_z) =j(\sigma\otimes \chi_z)^{-1}\in \mathbb{R}_{\geq 0}$. 
\end{proof}

Using the result above, we can now show that for every zero (resp. pole) of $\mu$, there must be another zero (resp. pole) of $\mu$, of the same multiplicity. 

\begin{lem}\label{lem:formmu}
    Let $\sigma$ be unitary. Then the function $\mu(\sigma,z)$ has the form
    \[
    \mu(\sigma,z)=c_{\mu}\prod_i (z-\lambda_i)^{n_i}(z^{-1}-\overline{\lambda_i})^{n_i},
    \]
    where $c_{\mu}\in\mathbb{R}_{>0}$, $\lambda_i\in\mathbb{C}$ and $n_i\in \mathbb{Z}$.
\end{lem}
\begin{proof}
By Lemma \ref{lem:unitcircle}, we have that $\mu(\sigma,z)\in\R_{\geq 0}$ for all $z$ on the unit circle. Since $\mu(\sigma,z)$ is a nonzero complex rational function in the variable $z$, we can write 
    \[
    \mu(\sigma,z)=c_{\mu}\prod_i(z-\lambda_i)^{n_i},
    \]
    where $c_{\mu}\in\mathbb{C}^{\times}$, $\lambda_i\in\mathbb{C}$, $n_i\in\mathbb{Z}$. For every negative power $n_i\in\mathbb{Z}_{<0}$ appearing in the factorization above, we may multiply $\mu(\sigma,z)$ by $(z-\lambda_i)^{-n_i}(z^{-1}-\overline{\lambda_i})^{-n_i}$, to obtain a function of the form
    \[
    \tilde{\mu}(\sigma,z)=c_{\mu}\prod_i(z-\lambda_i)^{n_i}\prod_j(z^{-1}-\gamma_j)^{\ell_j},
    \]
    where $\lambda_i,\gamma_j\in\mathbb{C}$ and $n_i,\ell_j\in\mathbb{Z}_{>0}$. Note that $\tilde{\mu}(\sigma,z)$ is again real and non-negative on the unit circle, and it suffices to show the statement of the lemma for $\tilde{\mu}(\sigma,z)$. For $\varphi\in \mathbb{R}$, we define $f(\varphi):=\tilde{\mu}(\sigma,e^{i\varphi})$. Then we have $f(\varphi)=\overline{f(\varphi)}$ for all $\varphi\in\mathbb{R}$, i.e.
    \[
    c_{\mu}\prod_i (e^{i\varphi}-\lambda_i)^{n_i} \prod_j (e^{-i\varphi}-\gamma_j)^{\ell_j}=\overline{c_{\mu}}\prod_i (e^{-i\varphi}-\overline{\lambda_i})^{n_i} \prod_j (e^{i\varphi}-\overline{\gamma_j})^{\ell_j}.
    \]
    Both sides in the equation are trigonometric polynomials, which have a unique factorization. So the set of $\overline{\lambda_i}$'s, where each $\overline{\lambda_i}$ appears $n_i$ times, equals the set of $\gamma_j$'s, where each $\gamma_j$ appears $\ell_j$ times. It follows that 
    \[
    \tilde{\mu}(\sigma,z)=c_{\mu}\prod_i(z-\lambda_i)^{n_i}(z^{-1}-\overline{\lambda_i})^{n_i}.
    \]
    Since $\mu$ is real and non-negative on the unit circle and not identically zero, we see that $c_{\mu}\in\mathbb{R}_{>0}$, which completes the proof.
\end{proof}

\subsection{An aspect of the Iwasawa decomposition}

In this subsection, we discuss an aspect of the Iwasawa decomposition. This will be necessary to study the limit of the intertwining operator $J_{\bar P | P}(\sigma\otimes\chi)$ and the function $\mu(\sigma,z)$ at $z=\infty$. We use the notations from Section \ref{sec:adjointness}. Recall the Iwasawa decomposition $G' = P' K' = K' P'$ 
from \cite[Theorem 5.3.4]{KaPr}. It lifts to the Iwasawa decomposition
\begin{equation}\label{eq:B.1}
G = P K = K P .
\end{equation}
As $P = U_P \rtimes M$, this means that for each $g \in G$ there exist $k_P(g) \in K, 
m_P(g) \in M$ and $u_P(g) \in U_P$ such that $g = u_P(g) m_P(g) k_P(g)$. The decomposition of $g$ along 
\eqref{eq:B.1} is unique up to $P \cap K$. Hence $k_P(g), m_P(g)$ and $u_P(g)$ are unique up to,
respectively, $P \cap K, M \cap K$ and $U_P \cap K$. 

For any subgroup $H$ of $G$ or $G$ we will write $H_x$ for the stabilizer of $x \in \mc B (G')$
in $H$. In particular $K' = G'_x$ and $K = G_x$. We choose a positive system $\Phi (G,S)^+$
in $\Phi (G,S)$, such that $\Phi (U_P,S) \subset \Phi (G,S)^+$. Let $U^+$ (respectively
$U^-$) be the subgroup of $G$ generated by the root subgroups $U_\beta$ with 
$\beta$ in $\Phi (G,S)^+$ (respectively in $\Phi (G,S)^- = -\Phi (G,S)^+$).

We note that the normalizer $N_G (S)$ equals $p_G^{-1}(N_{G'}(S'))$.

\begin{thm}\label{thm:B.1}
The product map
\[
\prod\nolimits_{\beta \in \Phi (G,S)^\pm} U_{\beta,x} \longrightarrow 
U^\pm_x = U^\pm \cap G_x
\]
is a homeomorphism. Further, $G_x$ can be obtained as
\[
G_x = U^-_x U^+_x N_G (S)_x = U^+_x U^-_x N_G (S)_x .
\]
\end{thm}
\begin{proof}
For $G'$, the theorem is the combination of \cite[Proposition 6.4.9]{BrTi} and
the proof of \cite[Lemma 7.7.3]{KaPr}. The theorem for $G$ follows from that for
$G'$ by the covering $p_G : G \to G'$ and the canonical lifting of each 
root subgroup of $G'$ to $G$.
\end{proof}

We will also work with the Lie algebra $\mr{Lie}(G')$, which can be
defined either as the $F$-points of Lie$(\mc G')$ or as the tangent space at $e$
of $G$ considered as manifold over the local field $F$. Since $p_G : G \to G'$ is 
a finite covering, we can identify Lie$(G')$ (defined using $F$-manifolds) with
$\mf{g} = \mr{Lie}(G)$. Similarly $\mr{Lie}(G'_x)$ can be viewed either as the 
$\mf o_F$-points of the Lie algebra of the $\mf o_F$-scheme $\mc G'_x$, or as the Lie 
algebra of $G'_x$ in the sense of manifolds over $\mf o_F$. Then $\mf g_x = \mr{Lie}(G_x)$ 
can be identified with Lie$(G'_x)$, and it is an $\mf o_F$-lattice in $\mf g$.
We write $\mf u_\beta = \mr{Lie}(U_\beta)$ and 
\[
\mf z = \mr{Lie}(Z_{G'}(S')) = \mr{Lie}(N_{G'}(S')) = \mr{Lie}(N_G (S)).
\]
By Lie$(H)_x$ we mean always Lie$(H_x)$. With these notations, 
Theorem \ref{thm:B.1} implies that 
\begin{align}\label{eq:B.2}
\mf g_x & = \bigoplus_{\beta \in \Phi (G,S)^-} \mf u_{\beta,x} \, \oplus
\bigoplus_{\beta \in \Phi (G,S)^+} \mf u_{\beta,x} \, \oplus \mf z_x \\
\nonumber & = \bigoplus_{\beta \in \Phi (U_{\bar P},S)} \mf u_{\beta,x} \, \oplus 
\Big( \bigoplus_{\beta \in \Phi (M,S)^-} \mf u_{\beta,x} \, \oplus 
\bigoplus_{\beta \in \Phi (M,S)^+} \mf u_{\beta,x} \, \oplus \mf z_x \, \Big) 
\oplus \bigoplus_{\beta \in \Phi (U_P,S)} \mf u_{\beta,x} \\
\nonumber & = \mr{Lie} (U_{\bar P,x}) \oplus \mr{Lie}(M_x) \oplus \mr{Lie}(U_{P,x}) .
\end{align}

\begin{prop}\label{prop:B.2}
Let $\bar u \in U_{\bar P}$ and write it as $\bar u = u m k$ 
with $u = u_P(\bar u) \in U_P, m = m_P(\bar u) \in M$ and $k = k_P(\bar u) \in G_x$.
Then $\delta_P (m_P(\bar u)) \leq 1$.
\end{prop}
\begin{proof}
Since $k \in G_x$, we have
\[
\bar u G_x \bar{u}^{-1} = G_{\bar{u} x} = G_{u m k x} = G_{u m x} = u G_{mx} u^{-1}.
\]
By \eqref{eq:B.2}, the Lie algebra of this group is
\begin{equation}\label{eq:B.3} 
\begin{split}
\mf{g}_{\bar u x} = \mr{Ad}(u) \mf{g}_{mx} & = \mr{Ad}(u) \big( \mr{Lie} (U_{\bar P,mx}) 
\oplus \mr{Lie}(M_{mx}) \oplus \mr{Lie}(U_{P,mx}) \big) \\
& = \mr{Ad}(u) \mr{Lie} (U_{\bar P,mx}) \oplus \mr{Ad}(u) \mr{Lie} (P_{mx}) .
\end{split}
\end{equation}
The decomposition $\mf g = \mr{Lie}(U_{\bar P}) \oplus \mf{p}$ gives a projection
$\mr{pr}_{\mf g / \mf p} : \mf g \to \mr{Lie}(U_{\bar P})$ with kernel 
$\mf p = \mr{Lie}(P)$. Clearly $\mr{pr}_{\mf g / \mf p} (\mf g_{\bar u x})$ contains
$\mr{Lie}(U_{\bar P,\bar u x}) = \mr{Ad}(\bar u) \mr{Lie}(U_{\bar P,x})$.
From \eqref{eq:B.3} we obtain, using $\mr{Ad}(u) \mr{Lie}(P_{mx}) \subset \mf{p}$:
\begin{equation}\label{eq:B.4}
\mr{Ad}(\bar u) \mr{Lie}(U_{\bar P,x}) \subset 
\mr{pr}_{\mf g / \mf p} (\mf g_{\bar u x}) = 
\mr{pr}_{\mf g / \mf p} \big( \mr{Ad}(u) \mr{Lie} (U_{\bar P,mx}) \big) .
\end{equation}
The Haar measures on the groups $G,M,U_P,P,U_{\bar P}, \bar P$ induce measures on their 
Lie algebras. By the unimodularity of the unipotent group $U_{\bar P}$, the lattice
$\mr{Ad}(\bar u) \mr{Lie}(U_{\bar P,x})$ in $\mr{Lie}(U_{\bar P})$ has the same volume
as $\mr{Lie}(U_{\bar P,x})$.

Since $u$ is unipotent, so is the $F$-linear operator Ad$(u) : \mf g \to \mf g$.
Moreover Ad$(u)$ stabilizes $\mf p$, so there exist a basis of $\mf p$ and a basis
of $\mr{Lie}(U_{\bar P})$ such that the matrix of Ad$(u)$ with respect to the
combined basis of $\mf g$ is unipotent and upper triangular. It follows that the
matrix of 
\begin{equation}\label{eq:B.5}
\mr{pr}_{\mf g / \mf p} \circ \mr{Ad}(u) : \mr{Lie}(U_{\bar P}) \to \mr{Lie}(U_{\bar P})
\end{equation}
with respect to the chosen basis is also unipotent and upper triangular. 
Therefore \eqref{eq:B.5} preserves volumes. Now we can interpret \eqref{eq:B.4}
in terms of volumes:
\begin{equation*}
\mr{vol} \big( \mr{Lie}(U_{\bar P,x}) \big) \leq 
\mr{vol} \big( \mr{Lie} (U_{\bar P,mx}) \big) .
\end{equation*}
This can be rewritten as 
\begin{equation}\label{eq:B.6}
\begin{aligned}
1 \leq \frac{\mr{vol} \big( \mr{Lie} (U_{\bar P,mx}) \big)}{\mr{vol} 
\big( \mr{Lie}(U_{\bar P,x}) \big)} & = \frac{\mr{vol} \big( \mr{Ad}(m) \mr{Lie} 
(U_{\bar P,x}) \big)}{\mr{vol} \big( \mr{Lie}(U_{\bar P,x}) \big)} \\
& = \big| \det \big( \mr{Ad}(m) : \mr{Lie}(U_{\bar P}) \to \mr{Lie}(U_{\bar P}) 
\big) \big|_F . 
\end{aligned}
\end{equation}
The modular function $\delta_{\bar P}$ of $\bar P$ equals $\delta_{\bar P'} 
\circ p_G$, because $\ker p_G$ is finite. By \cite[Lemma 1.2.1.1]{Sil}, we have
\begin{equation}\label{eq:B.7}
\delta_{\bar P'}(m') = \big| \det \big( \mr{Ad}(m') : \mr{Lie}(U_{\bar P}) \to 
\mr{Lie}(U_{\bar P}) \big) \big|_F \quad \text{for all } m' \in M'.
\end{equation}
The adjoint representation of $\bar P$ factors through $p_G : \bar P \to \bar P'$,
so the inequality \eqref{eq:B.6} says that $\delta_{\bar P}(m) \geq 1$.
Finally, we use that $\delta_{\bar P} = \delta_P^{-1}$ \cite[Lemme V.5.4]{Ren}.
\end{proof}

One may reformulate the conclusion of Proposition \ref{prop:B.2} as: 
conjugation by $m_P(\bar u)$ contracts $U_P$ and inflates $U_{\bar P}$. 
Thus the image of $\bar U$ in $M$ under the Iwasawa decomposition (well-defined 
up to $M \cap G_x$) sits inside $M$ in asymmetric way, 
like a semigroup but certainly not like a group. 

\begin{lem}\label{lem:B.3}
In the setting of Proposition \ref{prop:B.2}, suppose that $P$ is a maximal
proper parabolic subgroup of $G$ and that $\delta_P (m_P(\bar u)) = 1$.
Then $m_P(\bar u) \in M^1$.
\end{lem}
\begin{proof}
By construction $M^1 = p_G^{-1}(M'^1)$. As the Iwasawa decomposition of $G$ is
lifted from that for $G'$, it suffices to prove the lemma with $G'$ instead
of $G$. To use that simplification without many extra primes, we assume in
the remainder of the proof that $G = G'$.

The derived group $G_\der = \mc G_\der (F)$ also has an Iwasawa decomposition,
namely 
\begin{equation}\label{eq:B.9}
G_\der = U_P (M \cap G_\der) G_{\der,x}. 
\end{equation}
This is a restriction of the Iwasawa decomposition $G = U_P M G_x$. As 
$\bar u \in G_\der$, we can take $k_P(\bar u), m_P(\bar u)$ and $u_P(\bar u)$ in 
$G_\der$ as well, using \eqref{eq:B.9}.

By assumption $P \cap G_\der = (M \cap G_\der) \ltimes U_P$ is a maximal
proper parabolic subgroup of the semisimple group $G_\der$. Then the algebraic
group $Z(\mc M \cap \mc G_\der)^\circ$ is a onedimensional $F$-split torus. It
follows that $(M \cap G_\der) / (M \cap G_\der)^1$ is a free abelian group
of rank one, so isomorphic to $\Z$. 

On the right hand side of \eqref{eq:B.7} we can replace $\bar P$ by
$P$ or by $P \cap G_\der$, both give the same expression. Hence the modular
functions $\delta_P$ and $\delta_{P \cap G_\der}$ coincide on $M \cap G_\der$.
As $\delta_{P \cap G_\der}$ is trivial on $(M \cap G_\der)^1$, it defines 
a character
\begin{equation}\label{eq:B.8}
\Z \cong (M \cap G_\der) / (M \cap G_\der)^1 \to \R_{>0}^\times .
\end{equation}
The restrictions of elements of $\Phi (U_P,S)$ to $Z(M \cap G_\der)$ are
positive powers of a single nontrivial character $\alpha$. Then \eqref{eq:B.7}
shows that $\delta_{P \cap G_\der} |_{Z(M \cap G_\der)}$ is a positive power of 
$|\alpha|$. In particular $\delta_{P \cap G_\der}$ is nontrivial, so regarded
as a character of $\Z$ via \eqref{eq:B.8} it is injective. 
The assumptions of the lemma and Proposition \ref{prop:B.2} entail that
\[
\delta_{P \cap G_\der}(m_P(\bar u)) = \delta_P (m_P(\bar u)) = 1. 
\]
Therefore $m_P(\bar u)$ lies in the kernel of
\eqref{eq:B.8}, which means that 
\[
m_P(\bar u) \in (M \cap G_\der)^1 \subset M^1.
\]
That proves the lemma for one choice of $m_P(\bar u)$ in the Iwasawa decomposition
of $G$. But $m_P(\bar u)$ is unique up to the compact group $M \cap G_x$, so
$m_P(\bar u) \in M^1$ for any other choice of $m_P(\bar u)$ as well.
\end{proof}

\subsection{The limit at infinity}
In this subsection, we study the limit of the intertwining operator $J_{\bar P| P}(\sigma\otimes\chi)$ and the function $\mu(\sigma,z)$ when $z=\chi(h_\alpha^{\vee})$ goes to infinity.

By (\ref{eq:B.7}), the modulus character $\delta_P$ only takes values in the integer powers of $q$. By Proposition \ref{prop:B.2}, $\delta_P(m_P(\bar u))\leq 1$ for all $\bar u\in U_{\bar P}$. Hence, $U_{\bar P}$ is the disjoint union of the subsets 
\[
U_{\bar P}(n):=\set{\bar u \in U_{\bar P}: \delta_P(m_P(\bar u))=q^{-n}},
\]
where $n$ runs over all non-negative integers. The set 
\[U_{\bar P}(0)=\set{\bar u \in U_{\bar P}: \delta_P(m_P(\bar u))=1}\] is compact by \cite[Lemme II.3.4]{Wal} and $m_P(U_{\bar P}(0))\subset M^1$ by Lemma \ref{lem:B.3}. 

\begin{lem}\label{lem:maatU(n)}
    There exist $C>0$, $R>0$ and $d\in\mathbb{N}$ such that for every $n\geq 1$, 
    \[
    \mr{vol}(U_{\bar P}(n))\leq Cq^{nR}n^d.
    \]
\end{lem}
\begin{proof}
    In \cite[Lemme II.4.1]{Wal}, this statement (with $R=1/2$) is proven for $\delta_0:=\delta_{P_0}$, where $P_0$ is a minimal parabolic subgroup of $G$. So let $P_0=M_0U_P$ with $M_0\subset M$ a minimal Levi subgroup. Consider the Iwasawa decomposition $G=P_0K$, and for every $\bar u\in U_{\bar P}$, write 
    \[
    \bar u= u_0(\bar u)m_0(\bar u)k_0(\bar u),
    \]
    where $u_0(\bar u)\in U_P$, $m_0(\bar u)\in M_0$ and $k_0(\bar u)\in K$. Define, for all $n\geq 1$, the subsets 
    \[
    U_{\bar P, 0}(n):=\set{\bar u\in U_{\bar P}: \delta_0(m_0(\bar u))=q^{-n}}.
    \]
    Then, by \cite[Lemme II.4.1]{Wal}, there exist $C_0>0$ and $d_0\in\mathbb{N}$ such that for all $n\geq 1$, 
    \[
    \mr{vol}(U_{\bar P,0}(n))\leq C_0q^{n/2}n^{d_0}.
    \]
    By \cite[Lemme II.3.4]{Wal} and Proposition \ref{prop:B.2}, there exist $C_1>0$, $r>0$, such that for all $\bar u\in U_{\bar P}$,
    \[
    C_1\delta_P(m_P(\bar u))^r\leq \delta_0(m_0(\bar u)) \leq 1.
    \]
    Suppose $\bar u\in U_{\bar P}(n)$, i.e. $\delta_P(m_P(\bar u))=q^{-n}$. Then by the inequalities above, we get 
    \[
    C_1q^{-nr}\leq \delta_0(m_0(\bar u)) \leq 1.
    \]
    Write $C_1=q^{-b}$ for some $b\in\mathbb{R}$. Let $N\in\mathbb{N}$ such that $N\geq r+\max (b,0)$. Then for all $n\geq 1$,
    \[
    q^{-nN}\leq q^{-nr-b}=C_1q^{-nr}\leq \delta_0(m_0(\bar u)) \leq 1.
    \]
    Hence, $\bar u \in U_{\bar P}(n)$ can only be contained in one of the sets 
    \[
    U_{\bar P, 0}(nN), \ U_{\bar P, 0}(nN-1), \ \ldots, \ U_{\bar P, 0}(0).
    \]
    Moreover, each of the sets above has volume bounded above by $C_0q^{nN/2}(nN)^{d_0}$. Therefore
    \[
    \mr{vol}(U_{\bar P}(n))\leq (nN+1)C_0q^{nN/2}(nN)^{d_0}\leq n(N+1)C_0q^{nN/2}(nN)^{d_0}=Cq^{nR}n^d,
    \]
    where $C=(N+1)C_0N^{d_0}$, $R=N/2$ and $d=d_0+1$. 
\end{proof}

For $\chi \in X_{\nr}(M)$, we define the variable $r(\chi)\in \mathbb{R}$ determined by
\[
|\chi||_{M\cap G_{\rm{der}}}=\delta_P^{r(\chi)}|_{M\cap G_{\rm{der}}}.
\]
When this variable goes to infinity, the operator $J_{\bar P| P}(\sigma\otimes\chi)$ converges to an operator which is independent of $\chi$, as we will show below.

\begin{prop}\label{prop:infinitelimitexists} For $f\in I_{P\cap K}^K (\sigma)$, $g\in K$, 
\[
\lim_{r(\chi)\rightarrow \infty}J_{\bar P|P} (\sigma\otimes\chi)f(g)=\int_{U_{\bar P}(0)} f(\bar u g)\td \bar u,
\]
and this integral is independent of $\chi$. In particular, $\lim_{r(\chi)\rightarrow \infty}J_{\bar P|P} (\sigma\otimes\chi)$ exists.
\end{prop}
\begin{proof}
    For $\chi\in X_{\nr}(M)$, $f\in I_P^G(\sigma\otimes\chi)\cong I_{P\cap K}^K(\sigma)$ and $g\in K$, we have that 
\[J_{\bar P|P} (\sigma\otimes\chi)f(g)=\int_{U_{\bar P}} f(\bar u g)\td \bar u=\int_{U_{\bar P}(0)} f(\bar u g)\td \bar u+\int_{\cup_{n\geq 1}U_{\bar P}(n)} f(\bar u g)\td \bar u.\]
Recall that for $\bar u\in U_{\bar P}$, using the Iwasawa decomposition $G=PK$ we write
    \[
    \bar u = u_P(\bar u)m_P(\bar u)k_P(\bar u),
    \]
    where $u_P(\bar u)\in U_P$, $m_P(\bar u)\in M$ and $k_P(\bar u)\in K$. Then 
    \[
    f(\bar u g)=\chi\delta_P^{1/2}(m_P(\bar u))\sigma(m_P(\bar u)) f(k_P(\bar u) g).
    \]
    For $\bar u\in U_{\bar P}(0)$, we have that $\delta_P(m_P(\bar u))=1$. Then by Lemma \ref{lem:B.3}, $m_P(\bar u)\in M^1$, so also $\chi(m_P(\bar u))=1$. Hence 
    \[
    \int_{U_{\bar P}(0)} f(\bar u g)\td \bar u=\int_{U_{\bar P}(0)} \sigma(m_P(\bar u))f(k_P(\bar u) g)\td \bar u,
    \]
    which does not depend on $\chi$, since $k_P(\bar u)g\in K$ and $f\in I_{P\cap K}^K(\sigma)$.
    
    Now, let us analyze the integral 
    \[
    \int_{\cup_{n\geq 1}U_{\bar P}(n)} f(\bar u g)\td \bar u,
    \]
     in particular we want to show that this integral converges to zero in the limit $r(\chi)\rightarrow \infty$. For this, we follow the proof of \cite[Th\'eor\`eme IV.1.1]{Wal}. This proof shows that we can reduce the problem to showing that the integral 
     \[
     \int_{\cup_{n\geq 1} U_{\bar P}(n)} |\chi\delta_P^{r+1/2}(m_P(\bar u))|\td \bar u=\int_{\cup_{n\geq 1} U_{\bar P}(n)} \delta_P^{r(\chi)+r+1/2}(m_P(\bar u))\td \bar u,
     \]
where $r\in \mathbb{R}$ is some fixed number, converges to zero as $r(\chi)\rightarrow \infty$. By Lemma \ref{lem:maatU(n)}, there exist $C>0$, $R>0$ and $d\in\mathbb{N}$ such that for each $n\geq 1$, 
\[
\mr{vol}(U_{\bar P}(n))\leq Cq^{nR}n^d.
\]
Hence, we obtain that 
\begin{align*}
\int_{\cup_{n\geq 1} U_{\bar P}(n)} \delta_P^{r(\chi)+r+1/2}(m_P(\bar u))\td \bar u&\leq \sum_{n\geq 1} q^{-n(r(\chi)+r+1/2)}Cq^{nR}n^d\\
&=C\sum_{n\geq 1} q^{-n(r(\chi)+r')}n^d, 
\end{align*}
where $r'\in\mathbb{R}$ is independent of $n$. This sum converges for $r(\chi)$ large enough, and it converges to zero when $r(\chi)\rightarrow \infty$.  
\end{proof}

 Recall that for every $\chi$, we may view $J_{\bar P | P}(\sigma\otimes\chi)$ as a linear operator 
    \[
    J_{\bar P | P}(\sigma\otimes\chi):I_{P\cap K}^K(V_\sigma)\rightarrow I_{\bar P\cap K}^K(V_\sigma),
    \]
    i.e. these are all linear operators on the same vector space, with the same range. If $z=\chi(h_\alpha^{\vee})\rightarrow \infty$, then $r(\chi)\rightarrow \infty$. So by Proposition \ref{prop:infinitelimitexists}, $\lim_{z\rightarrow\infty}J_{\bar P | P}(\sigma\otimes\chi)$ exists and is given by 
    \[
    \lim_{z\rightarrow\infty}J_{\bar P | P}(\sigma\otimes\chi)=\int_{U_{\bar P}(0)}\lambda(\bar u)^{-1} \td \bar u:I_{P\cap K}^K(V_\sigma)\rightarrow I_{\bar P\cap K}^K(V_\sigma).
    \]
    One can also directly prove that the image of this operator lies in $I_{\bar P\cap K}^K(V_\sigma)$, as we will do below. 
    \begin{lem}
        Let $f\in I_{P\cap K}^K(V_\sigma)$ and let $\bar f:=\int_{U_{\bar P}(0)}\lambda(\bar u)^{-1}f \td \bar u$. Then $\bar f\in I_{\bar P\cap K}^K(V_\sigma)$. 
    \end{lem}
    \begin{proof}
    Let $f\in I_{P\cap K}^K(V_\sigma)$, so $f:K\rightarrow V_\sigma$ is a function right invariant under some compact open subgroup and $f(mug)=\sigma(m)f(g)$ for all $g\in K$, $m\in M\cap K$, $u\in U_P\cap K$. We want to show that for $m\in M\cap K$, $\bar u \in U_{\bar P}\cap K$, $g\in K$, we have $\bar f(m\bar u g)=\sigma(m)\bar f(g)$. First, let us show that $U_{\bar P} \cap K\subset U_{\bar P}(0)$. For $\bar u \in U_{\bar P}\cap K$, write 
    \[
    \bar u = u_P(\bar u)m_P(\bar u ) k_P(\bar u )
    \]
    using the Iwasawa decomposition. We may choose $k_P(\bar u)=\bar u \in K$. Then $m_P(\bar u)=1$, hence $\delta_P(m_P(\bar u))=1$ so $\bar u\in U_{\bar P}(0)$. 

    Now let $m\in M\cap K$, $\bar u \in U_{\bar P}\cap K$, $g\in K$. Then 
    \begin{align*}
    \bar f(m\bar u g)&=\int_{U_{\bar P}(0)}f(\bar u_1m\bar ug)\td \bar u_1=\int_{U_{\bar P}(0)}f(mm^{-1}\bar u_1m\bar u g)\td\bar u_1\\
    &=\sigma(m)\int_{U_{\bar P}(0)}f(m^{-1}\bar u_1m\bar ug)\td\bar u_1.
    \end{align*}
    We have that 
    \[
    m^{-1}\bar u_1m= (m^{-1}u_P(\bar u_1)m)(m^{-1}m_P(\bar u_1)m)(m^{-1}k_P(\bar u_1)m),
    \]
    therefore
    \[
    \delta_P(m_P(m^{-1}\bar u_1 m))=\delta_P(m^{-1}m_P(\bar u_1)m)=1,
    \]
    since $m\in K$ and $\bar u_1\in U_{\bar P}(0)$. So $m^{-1}\bar u_1 m\in U_{\bar P}(0)$ and this gives
    \[
    \Ad(m^{-1}):U_{\bar P}(0)\rightarrow U_{\bar P}(0),
    \]
    which preserves measures since $m\in K\cap M$, so $\delta_P(m)=1$. Thus
    \begin{align*}
        \bar f(m\bar u g)&=\sigma(m)\int_{U_{\bar P}(0)}f(m^{-1}\bar u_1m\bar ug)\td\bar u_1=\sigma(m)\int_{U_{\bar P}(0)}f(\bar u_1\bar ug)\td\bar u_1.
    \end{align*}
    We can write
    \[
    \bar u_1\bar u = u_P(\bar u_1)m_P(\bar u_1)(k_P(\bar u_1)\bar u),
    \]
    so $\bar u_1\bar u\in U_{\bar P}(0)$, and we may change the variable $\bar u_1 \bar u$ to $\bar u_2\in U_{\bar P} (0)$, by right invariance of the Haar measure on $U_{\bar P}$. So we obtain 
    \[
    \bar f(m\bar u g)= \sigma(m)\int_{U_{\bar P}(0)}f(\bar u_2g)\td\bar u_2=\sigma(m)\bar f (g),
    \]
    hence $\bar f\in I_{\bar P\cap K}^K(V_\sigma)$.
    \end{proof}

Since 
    \[
    J_{P|\bar P}(\sigma\otimes \chi)\circ J_{\bar P | P}(\sigma\otimes\chi)=j(\sigma\otimes\chi)=\mu(\sigma\otimes\chi)^{-1}
    \]
    is a rational function in the variable $z=\chi(h_\alpha^\vee)$, it has only finitely many poles and zeros. So for $z=\chi(h_\alpha^\vee)$ large enough, $J_{\bar P| P}(\sigma\otimes\chi)$ is invertible. We can use this to show that $\lim_{z\rightarrow\infty}J_{\bar P| P}(\sigma\otimes\chi)$ is invertible. 

\begin{lem}\label{lem:infinitelimisinvertible}
    The operator
    \[
    A:=\int_{U_{\bar P}(0)}\lambda(\bar u)^{-1}\td \bar u : I_{P\cap K}^K(V_\sigma) \rightarrow I_{ \bar P \cap K}^K(V_\sigma)
    \]
    is invertible.
\end{lem}
\begin{proof}
    It suffices to show that for every compact open subgroup $C$ of $K$, the operator 
    \[
    A^C:=\int_{U_{\bar P}(0)}\lambda(\bar u)^{-1}\td \bar u : I_{P\cap K}^K(V_\sigma)^C \rightarrow I_{\bar P \cap K}^K(V_\sigma)^C,
    \]
     which is a map of finite-dimensional vector spaces, is bijective. Let $z=\chi(h_\alpha^\vee)$ be large enough such that $J_{\bar P | P}(\sigma\otimes\chi)$ is invertible. Let $f\in I_{P\cap K}^K (V_\sigma)^C$, but view it as an element of $I_P^G(\sigma\otimes\chi)$. Let $g\in K$. Then
    \[
    \int_{U_{\bar P}}\left[\int_{U_{\bar P}(0)}f(\bar u_1 \bar u_2g)\td \bar u_1\right]\td \bar u_2=\mr{vol}(U_{\bar P}(0))\int_{U_{\bar P}}f(\bar u g)\td \bar u,
    \]
    hence 
    \[
    \int_{U_{\bar P}}\lambda(\bar u)^{-1}\td\bar u \circ \int_{U_{\bar P}(0)}\lambda(\bar u)^{-1}\td \bar u = \mr{vol}(U_{\bar P}(0))J_{\bar P | P}(\sigma\otimes \chi).
    \]
    Since $J_{\bar P | P}(\sigma\otimes \chi)$ is invertible, $\int_{U_{\bar P}(0)}\lambda(\bar u)^{-1}\td \bar u$ is injective, hence $A^C$ is bijective since it is a map of finite dimensional vector spaces. We can do this for every $C$, so $A$ is invertible.
\end{proof}

We can use the results above to show that $\mu(\sigma,z)=\mu(\sigma\otimes\chi)=j(\sigma\otimes\chi)^{-1}$ is not zero at $z=\infty$ and $z=0$.

\begin{prop}\label{prop:muatinfinity}
The function $\mu(\sigma,z)$ is not zero at $z=\infty$ and $z=0$, i.e. 
\[\lim_{z\rightarrow\infty} j(\sigma\otimes\chi)\not=\infty \textnormal{ and } \lim_{z\rightarrow 0} j(\sigma\otimes\chi)\not=\infty.\]
\end{prop}
\begin{proof}
We will first show that $\lim_{z\rightarrow\infty} j(\sigma\otimes\chi)\not=\infty$. Let $f_{\chi}\in I_{\bar P}^G(\sigma\otimes\chi)$ be determined by $\supp(f_\chi)=\bar P U_{P,r}$, $f_\chi(1)=v$, and $f_\chi$ is right $U_{P,r}$-invariant. Then
    \[
    J_{P|\bar P}(\sigma\otimes\chi)f_\chi(1)=\int_{U_P}f_\chi(u)\td u=\int_{U_{P,r}}f_\chi(u)\td u=\mr{vol}(U_{P,r})v,
    \]
    which is independent of $\chi$. Let $z=\chi(h_\alpha^\vee)$ be large enough such that $J_{\bar P | P}(\sigma\otimes\chi)$ is invertible, and write $g_\chi=J_{\bar P | P}(\sigma\otimes \chi)^{-1}f_\chi\in I_P^G(\sigma\otimes\chi)$. Then
    \[
    j(\sigma\otimes\chi)g_{\chi}(1)=J_{P|\bar P}(\sigma\otimes\chi)\circ J_{\bar P|P}(\sigma\otimes\chi)g_\chi(1)=J_{P| \bar P}(\sigma\otimes\chi)f_\chi(1)=\mr{vol}(U_{P,r})v. 
    \]
    By Proposition \ref{prop:infinitelimitexists} and Lemma \ref{lem:infinitelimisinvertible}, $\lim_{z\rightarrow\infty}J_{\bar P | P}(\sigma\otimes\chi)$ exists and is invertible. Hence, $\lim_{z\rightarrow\infty}g_\chi=g\in I_{P\cap K}^K(V_\sigma)$ exists. 
    Hence, 
    \[
    \lim_{z\rightarrow\infty}j(\sigma\otimes\chi)g_{\chi}(1)=\lim_{z\rightarrow\infty}j(\sigma\otimes\chi)g(1)=\mr{vol}(U_{P,r})v,
    \]
    so $\lim_{z\rightarrow\infty}j(\sigma\otimes\chi)\not=\infty$. 

    So far, we showed that $\mu(\sigma,\infty)\not=0$.
    By Lemma \ref{lem:supercuspidalunitary}, $\sigma=\tilde\sigma\otimes\tilde\chi$ for some unitary $\tilde\sigma\in \Irr(M)_{[M,\sigma]}$ and $\tilde\chi\in\Hom(M,\mathbb{R}_{>0})$. By (\ref{eq:mutau}), $\mu(\sigma,0)=\mu(\tilde\sigma,0)$ and $\mu(\sigma,\infty)=\mu(\tilde\sigma,\infty)$.
    By Lemma \ref{lem:formmu}, $\mu(\tilde\sigma,\infty)\not=0$ implies that $\mu(\tilde\sigma,0)\not =0$, i.e. $\mu(\sigma,0)\not=0$.
\end{proof}

\subsection{Zeros of the \texorpdfstring{$\mu$}{mu}-function}
\label{par:zeros}

In this subsection, we analyze the zeros of the $\mu$-function. Since $M$ is maximal, we have $[N_G(M):M]\leq 2$, and we will distinguish between two cases. 

Recall that Bernstein's geometric lemma \cite[Th\'eor\`eme VI.5.1]{Ren} provides 
a filtration of $J^G_{\bar P} I_P^G (\sigma)$. The associated subquotients have the form 
\begin{equation}\label{eq:4.1}
\mc F_{w,\bar P} / \mc F_{w+,\bar P} = I^M_{M \cap w P w^{-1}}
(w \cdot J^M_{M \cap w^{-1} \bar P w} \sigma ) 
\end{equation}
where $w \in W(M,S) \backslash W(G,S) / W(M,S) \cong \bar P \backslash G / P$.
By the cuspidality of $\sigma$, such a subquotient vanishes whenever (a representative
of) $w$ does not normalize $M$, for then $M \cap w^{-1} \bar P w$ is a proper parabolic
subgroup of $M$. As $N_{W(G,S)}(W(M,S)) / W(M,S) \cong N_G (M) / M$, Bernstein's
geometric lemma reduces to
\begin{align}\label{eq:4.2}
0 \to s_\alpha \cdot \sigma \to & J^G_{\bar P} I_P^G (\sigma) \to \sigma \to 0
\qquad N_G (M) / M = \{1,s_\alpha\} ,\\
\label{eq:4.3} & J^G_{\bar P} I_P^G (\sigma) \cong \sigma \hspace{18mm} N_G (M) = M .
\end{align}
In the case where $N_G(M)/M$ is trivial, it turns out that $\mu$ is constant. In the case where $N_G(M)/M$ is not trivial, we will prove that $\mu$ has at most two zeros in $\mathbb{C}^\times$.

\begin{lem}\label{lem:munozero}
    If $\mu(\sigma,z)$ has no zeros on $\mathbb{C}^\times$, then $\mu(\sigma,z)$ is constant. 
\end{lem}
\begin{proof}
    Take a unitary $\tilde\sigma\in \Irr(M)_{[M,\sigma]}$. Recall from Lemma \ref{lem:formmu} that we can write 
    \[
    \mu(\tilde\sigma,z)=c_{\mu}\prod_i (z-\lambda_i)^{n_i}(z^{-1}-\overline{\lambda_i})^{n_i}.
    \]
    Suppose $\mu(\tilde\sigma,z)$ has no zeros in $\mathbb{C}^\times$, then there are no factors with $n_i$ positive. So there can only be factors with $n_i$ negative. However, by Proposition \ref{prop:muatinfinity}, $\mu(\tilde\sigma,\infty)\not = 0$, so $\mu(\tilde\sigma,z)$ is constant, hence $\mu(\sigma,z)$ is also constant by (\ref{eq:mutau}).
\end{proof}

\begin{lem}\label{lem:muconstant}
If $N_G(M)/M$ is trivial, then $\mu(\sigma,z)$ is constant.
\end{lem}
\begin{proof}
We follow the proof of \cite[Th\'eor\`eme IV.1.1]{Wal}, where it is shown that 
$J_{\bar P | P}(\sigma \otimes \chi)$ depends rationally on $\chi \in X_\nr (M)$.
Since we are working with the covering $p : G \to G'$, we define $A_M$ as $Z(M) \cap S$.
This is not an algebraic torus like in \cite{Wal}, but $p(A_M)$ has finite index in
$Z(M') \cap S'$, which is good enough to make \cite[\S 4.1]{Wal} work. By \eqref{eq:4.1}
and \eqref{eq:4.3}, the representation $\mc F_{w,\bar P} / \mc F_{w+,\bar P}$ is nonzero
if and only if $w = 1$. Moreover $A_M$ acts by a unique character on the irreducible
representation $\sigma$. 

Then \cite[p. 281]{Wal} simplifies considerably: we have $\mathcal{E}xp_w = \emptyset$
if $w \neq 1$, $|\mathcal{E}xp_1| = 1$ and $R = b = 1$. Next the definitions of $j$ and
$J$ in \cite[p. 281--282]{Wal} show that via Frobenius reciprocity
\[
J_{\bar P | P}(\sigma \otimes \chi) : I_P^G (\sigma \otimes \chi) \to 
I_{\bar P}^G (\sigma \otimes \chi)
\]
corresponds to the isomorphism $J^G_{\bar P} I_P^G (\sigma \otimes \chi) \cong 
\sigma \otimes \chi$ from \eqref{eq:4.3}. In particular $J_{\bar P | P}(\sigma 
\otimes \chi)$ is nonzero and regular for all $\chi \in X_\nr (M)$. The same holds
for $J_{P | \bar P}(\sigma \otimes \chi)$. Therefore $\mu (\sigma,z) = \mu 
(\sigma \otimes \chi)$ is also nonzero and regular. Then Lemma \ref{lem:munozero}
shows that it is constant.
\end{proof}

So from now on, let us assume that $[N_G(M):M]=2$ and that $\mu$ has a zero, otherwise there is nothing left to prove. Let $s_\alpha$ be the unique nontrivial element of $N_G(M)/M$, with representative $n_\alpha \in N_G(M)$. Then $s_\alpha$ acts on $\Irr(M)$ via $(s_\alpha\cdot \pi )(m)=\pi(n_\alpha^{-1}m n_\alpha)$. 

\begin{prop}\label{prop:zeros}
    Suppose that $\mu(\sigma,z)$ is not constant. Then there exists a unitary $\sigma_0\in \Irr(M)_{[M,\sigma]}$ such that the following statements hold:
    \begin{enumerate}
        \item[(i)] $\mu(\sigma_0)=0$;
        \item[(ii)] $\mu(\sigma_0,z)=0$ implies that $z=1$ or $z=-1$;
        \item[(iii)] if $\mu(\sigma_0,z)=0$, the zero at $z$ is of order $2$.
    \end{enumerate}
\end{prop}
\begin{proof}
Since $\mu(\sigma,z)$ is not constant, $\mu(\tau)=0$ for some $\tau\in \Irr(M)_{[M,\sigma]}$,
by Lemma \ref{lem:munozero}. We claim that $s_\alpha \cdot \tau$ is isomorphic to 
$\tau$.

Suppose, on the contrary, that $s_\alpha \tau \not\cong \tau$.
As $\tau$ and $s_\alpha \tau$ are both supercuspidal, 
\eqref{eq:4.2} splits \cite[Lemme VI.3.6]{Ren} and becomes
\begin{equation}\label{eq:4.4}
J^G_{\bar P} I_P^G (\tau) \cong \mc F_{s_\alpha,\bar P} / \mc F_{s_\alpha +,\bar P} \oplus
\mc F_{1,\bar P} / \mc F_{1+,\bar P} \cong s_\alpha \tau \oplus \tau .
\end{equation}
Let $\C (X_\nr (M))$ be the field of rational functions on the complex torus $X_\nr (M)$.
Let $A_M$ be as in the proof of Lemma \ref{lem:muconstant} and let $\chi_\tau$ be the
$A_M$-character of $\tau$. We define $\mu_B : A_M \to \C (X_\nr (M))^\times$ as the canonical
homomorphism given by $\mu_B (a) (\chi) = \chi (a)$. In the notations from \cite[p. 280]{Wal},
we have 
\[
\mathcal{E}xp (\tau) = \{\chi_\tau\},\; \mathcal{E}xp_1 = \{\chi_\tau \mu_B \},\;
\mathcal{E}xp_{s_\alpha} = \{ s_\alpha (\chi_\tau) s_\alpha (\mu_B)\} 
\] 
and $\mathcal{E}xp_w$ is empty for all other $w$.
In \cite[Corollaire IV.1.2]{Wal} an upper bound $D$ for the pole order of $J_{\bar P | P}
(\tau \otimes \chi)$ at $\chi = 1$ is obtained. If $s_\alpha (\chi_\tau) \neq \chi_\tau$,
then $D = 0$ because it is a sum over an empty index set. If $s_\alpha (\chi_\tau) =
\chi_\tau$, then $D$ equals $d(\chi_\tau,s_\alpha)$, which is defined on \cite[p. 281]{Wal}.
In view of \eqref{eq:4.4}, $d(\chi_\tau,s_\alpha)$ is an integer such that, for all
$a \in A_M$, 
\[
\big( (s_\alpha \tau \otimes s_\alpha \mu_B)(a) - (\chi_\tau \otimes s_\alpha \mu_B) (a) 
\big)^{d(\chi_\tau,s_\alpha)} \text{ annihilates } s_\alpha \tau \otimes \C (X_\nr (M)) .
\]
But $s_\alpha \tau (a) = \chi_\tau (a) \mr{id}$ for all $a \in A_M$, so we may take
$d(\chi_\tau,s_\alpha) = 0$. 

So we have $D=0$ irrespective of $s_\alpha (\chi_\tau)$, and \cite[Corollaire IV.1.2]{Wal}
says that $J_{\bar P | P}(\tau \otimes \chi)$ is regular at $\chi = 1$. 
The same goes for $J_{P | \bar P}(\tau \otimes \chi)$. 
Hence $j(\tau \otimes \chi)$ is regular at $\chi = 1$, and $\mu (\tau) \neq 0$.
This contradicts the definition of $\tau$. We conclude that $s_\alpha \tau$ must be
isomorphic to $\tau$.

By Lemma \ref{lem:supercuspidalunitary}, there is a bijection
\[
\Irr(M)^\unit_{[M,\sigma]}\times \Hom(M,\mathbb{R}_{>0})\rightarrow \Irr(M)_{[M,\sigma]}: (\pi,\chi)\mapsto \pi\otimes \chi.
\]
This bijection is $s_\alpha$-equivariant by construction,
\[s_\alpha\cdot(\pi,\chi)=(s_\alpha\pi,s_\alpha\chi)\mapsto s_\alpha \pi\otimes s_\alpha\chi=s_\alpha\cdot(\pi\otimes \chi).\]
So $\tau=\sigma_0\otimes \chi_0$ for a unitary representation $\sigma_0$ and $\chi_0\in\Hom(M,\mathbb{R}_{>0})$, with $s_\alpha\sigma_0\cong \sigma_0$ and $s_\alpha\chi_0 \cong \chi_0$.

We will first prove \textit{(ii)} and use that result to prove \textit{(i)}. For $\mu(\sigma_0\otimes\chi)=\mu(\sigma_0,z)$ to be zero, we must have $s_\alpha (\sigma_0\otimes \chi)\cong \sigma_0\otimes \chi$ \cite[283]{Wal}. We have $s_\alpha \sigma_0\cong \sigma_0$, so $s_\alpha(\sigma_0\otimes\chi)\cong \sigma_0\otimes s_\alpha \chi$. Therefore, $\mu(\sigma_0\otimes\chi)=0$ implies that $\sigma_0\cong \sigma_0\otimes \chi (s_\alpha\chi)^{-1}$, so $\chi(s_\alpha\chi)^{-1}(h_\alpha^\vee)=1$. Write $\chi(h_\alpha^\vee)=z$, then $(s_\alpha\chi)(h_\alpha^\vee)=z^{-1}$, so $\chi(s_\alpha\chi)^{-1}(h_\alpha^\vee)=z^2$. Hence, $\mu(\sigma_0\otimes\chi)=\mu(\sigma_0,z)=0$ implies that $z^2=1$, hence $z=\pm 1$. This proves \textit{(ii)}. 

Recall that $\mu(\tau)=\mu(\sigma_0\otimes\chi_0)=0$, where $\chi_0\in\Hom(M,\mathbb{R}_{>0})$. By the above, we obtain that $\chi_0(h_\alpha^\vee)=1$, hence $\mu(\sigma_0)=0$. This proves \textit{(i)}. 

It remains to prove \textit{(iii)}. Suppose $\mu(\sigma_0,z)=0$, then $z=\pm1$ by \textit{(ii)}. In particular $z\in S^1$, so by the proof of Lemma \ref{lem:unitcircle}, there exists a unitary character $\chi_z\in X_\nr(M)$ such that $\chi_z(h_\alpha^\vee)=z$. Then $\sigma_0\otimes\chi_z$ is unitary with $\mu(\sigma_0\otimes\chi_z)=0$. This means that at least one of $J_{\bar P|P}(\sigma_0\otimes\chi_z)$ and $J_{P |\bar P }(\sigma_0\otimes\chi_z)$ is singular. Since $\sigma_0\otimes\chi_z$ is unitary (in particular Hermitian), we have the adjointness relation from Theorem \ref{thm:A.5}.a with $\overline{\sigma_0\otimes\chi_z}^\vee$ replaced by $\sigma_0\otimes\chi_z$. So
 \[
 \langle J_{\bar P|P}(\sigma_0\otimes \chi_z) f_1, f_2\rangle =\langle f_1, J_{P|\bar P}(\sigma_0\otimes \chi_z)f_2\rangle
 \]
 for all $f_1\in I_P^G(\sigma_0\otimes\chi_z)$, $f_2\in I_{\bar{P}}^G(\sigma_0\otimes \chi_z)$. Thus, singularity of one of the two intertwining operators $J_{\bar P|P}$ and $J_{P|\bar P}$ at $\sigma_0\otimes\chi_z$ automatically implies singularity of the other one at $\sigma_0\otimes\chi_z$. By \cite[Corollaire IV.1.2]{Wal}, both singularities have order $1$. Since 
 \[
 \mu(\sigma_0\otimes\chi)^{-1}=j(\sigma_0\otimes\chi)=J_{P|\bar P}(\sigma_0\otimes\chi)\otimes J_{\bar P|P}(\sigma_0\otimes\chi),
 \]
we obtain that the zero of $\mu$ at $\sigma_0\otimes\chi_z$ has order $2$. This proves \textit{(iii)}.
\end{proof}

\subsection{Hermitian techniques to find the poles}

In this subsection, we locate the poles of $\mu(\sigma,z)$ using techniques with Hermitian and unitary representations. We will prove the main result of this paper, which is a formula for $\mu(\sigma,z)$ in terms of poles and zeros. Recall that in Lemma \ref{lem:formmu}, we have already written the function $\mu(\sigma,z)$ in a form that is closer to the desired result. We also showed in Lemma \ref{lem:munozero} that $\mu(\sigma,z)$ is constant if it has no zeros, which is satisfied when $[N_G(M):M]=1$, see Lemma \ref{lem:muconstant}. 

Therefore, from now on, we assume that $[N_G(M):M]=2$ and that $\mu(\sigma,z)$ is not constant. So $N_G(M)/M$ has a unique nontrivial element $s_\alpha$. For simplicity, we assume in this section that $\sigma=\sigma_0$, where $\sigma_0$ is as in Proposition \ref{prop:zeros}. So $\sigma$ is unitary with $s_\alpha\sigma\cong \sigma$ and $\mu(\sigma)=0$. There is possibly another zero of $\mu(\sigma,z)$ at $z=-1$ and each zero is of order $2$. For each zero of $\mu$, we will find two poles in $\mathbb{R}^\times$. To do this, it suffices to find one pole for each zero, because by Lemma \ref{lem:formmu}, one obtains the other pole for free. We already showed in Proposition \ref{prop:muatinfinity} that $\mu(\sigma,z)$ is not zero at $z=\infty$ and $z=0$, so we can conclude that $\mu(\sigma,z)$ does not have any more poles in $\mathbb{C}^\times$.

We will use arguments concerning Hermitian and unitary representations to locate the poles. See \cite[\S IV]{Ren} for definitions and useful results about Hermitian and unitary representations. Let $\chi\in\Hom(M,\R_{>0})\subset X_{\nr}(M)$, then $\overline{\chi}=\chi$ and $\overline{\chi}^\vee=\chi^\vee=\chi^{-1}$. By unitarity of $\sigma$, we have $\overline{\sigma\otimes\chi}^\vee\cong\sigma\otimes\chi^{-1}$. So $\sigma\otimes\chi$ is only Hermitian if $\sigma\otimes\chi\cong\sigma\otimes\chi^{-1}$, or equivalently if $\chi^2\in X_{\nr}(M,\sigma)$. We will show that the induced representation $I_P^G(\sigma\otimes\chi)$ is Hermitian whenever $\chi\in\Hom(M/Z(G),\mathbb{R}_{>0})$ and $\mu(\sigma\otimes\chi)\not=\infty$. This will allow us to construct a continuous family of Hermitian forms on a family of Hermitian representations, which will help us to show that somewhere in this family, a pole of $\mu(\sigma\otimes\chi)$ must occur.

\begin{lem} \label{lem:renardirredsubquotients}
	Let $L$ be a Levi subgroup of $G$ and $\tau\in\Irr(L)$, $w\in N_G(L)$. Let $Q=LV$ be a parabolic subgroup of $G$. Then $I_Q^G(\tau)$ and $I_Q^G(w\cdot \tau)$ have the same irreducible subquotients with the same multiplicities. In particular, if $I_Q^G(\tau)$ is irreducible, then $I_Q^G(\tau)\cong I_Q^G(w\cdot \tau)$.
\end{lem}
\begin{proof}
	This follows from \cite[Theorem 2.9, Remark 2.10]{BerZelInducedRep}. See also \cite[Th\'eor\`eme VI.5.4]{Ren}.
\end{proof}

\begin{lem}\label{lem:inducedrepishermitian}
  Suppose that $s_\alpha(\chi)=\chi^{-1}$ and that $I_P^G(\sigma\otimes\chi)$ is irreducible. Then $I_P^G(\sigma\otimes\chi)$ is Hermitian.
\end{lem}
\begin{proof}
    Since $s_\alpha(\sigma\otimes\chi)\cong \sigma\otimes\chi^{-1}$, we have by Lemma \ref{lem:renardirredsubquotients} that $I_P^G(\sigma\otimes\chi)\cong I_P^G(\sigma\otimes\chi^{-1})$. Taking the complex conjugate representation commutes with parabolic induction. By \cite[\S VI.1.2]{Ren}, taking the contragredient representation also commutes with parabolic induction. So we get
    \[
    I_P^G(\sigma\otimes\chi)\cong I_P^G(\sigma\otimes\chi^{-1})\cong I_P^G(\overline{\sigma\otimes\chi}^\vee)\cong \overline{I_P^G(\sigma\otimes\chi)}^\vee,    
    \]
    where the middle isomorphism follows from unitarity of $\sigma$, so $\overline{\sigma\otimes\chi}^\vee\cong\sigma\otimes\chi^{-1}$. This completes the proof that $I_P^G(\sigma\otimes\chi)$ is Hermitian.
\end{proof}

The condition $s_\alpha\chi=\chi^{-1}$ from Lemma \ref{lem:inducedrepishermitian} is satisfied by unramified characters $\chi\in \Hom(M/Z(G),\mathbb{R}_{>0})$, as shown below.

\begin{lem}
    If $\chi\in \Hom(M/Z(G),\mathbb{R}_{>0})$, then $s_\alpha\chi=\chi^{-1}$.
\end{lem}
\begin{proof}
    The unramified character $\chi$ factors through $M/M^1$, so 
    \[\chi:M/M^1Z(G)\rightarrow\mathbb{R}_{>0}.\] Because $M$ is a maximal Levi subgroup, we have 
    \[
    M/M^1\cong \mathbb{Z}^r \textnormal{ and } Z(G)/Z(G)_\cpt\cong \mathbb{Z}^{r-1}.
    \] 
    This is because the rank $r$ of $M/M^1$ equals the dimension of the maximal split torus in $M$, which is one higher than the dimension of the maximal split torus in $G$, i.e. the rank of $Z(G)/Z(G)_\cpt$, which is the same as the rank $r-1$ of $G/G^1$. Thus, we obtain the quotient
    \[M/M^1Z(G)=(M/M^1)/(Z(G)/Z(G)_\cpt)\cong \mathbb{Z}\times H,\] with $H$ a finite group. So we can consider $\chi:\mathbb{Z}\times H\rightarrow \mathbb{R}_{>0}$, which must be trivial on $H$ since it is a finite group. So $\chi:\mathbb{Z}\rightarrow \mathbb{R}_{>0}$ is of the form $n\mapsto x^n$ where $x\in\mathbb{R}_{>0}$. Now $s_\alpha$ acts on $\mathbb{Z}$ as a nontrivial group homomorphism and an involution, so the only possibility is multiplication by $-1$. This induces the action on $\chi$, which we conclude is given by $s_\alpha\chi=\chi^{-1}$.
\end{proof}

Under the assumptions of Lemma \ref{lem:inducedrepishermitian}, we will now construct a Hermitian form on $I_P^G(\sigma\otimes\chi)$, using the intertwining operator $J_{\bar{P}|P}(\sigma\otimes\chi)$. Write $\lambda(s_\alpha)$ for translation by $s_\alpha$. Fix an isomorphism of $M$-representations $\varphi_\sigma:\sigma\rightarrow s_\alpha\cdot \sigma$. Moreover, let 
\[
J'_{\bar{P}|P}(\sigma\otimes\chi):=(\chi(h_\alpha^\vee)-1)J_{\bar{P}|P}(\sigma\otimes\chi), \ J'_{P|\bar{P}}(\sigma\otimes\chi):=(\chi^{-1}(h_\alpha^\vee)-1)J_{P|\bar{P}}(\sigma\otimes\chi).
\]
Then both $J'_{\bar{P}|P}(\sigma\otimes\chi)$ and $J'_{P|\bar{P}}(\sigma\otimes \chi)$ do not have a pole at $\chi=1$. They are invertible for $\chi$ in a neighbourhood of $\chi=1$ in $X_{\nr}(M)$. They are each other's inverse, up to a scalar which depends rationally on $\chi$.
Consider, for $\chi\in \Hom(M/Z(G),\mathbb{R}_{>0})$ the composition
\begin{align*}
&I_P^G(\sigma\otimes\chi)\xrightarrow{J'_{\bar{P}|P}(\sigma\otimes\chi)} I_{\bar{P}}^G(\sigma\otimes\chi)\xrightarrow{\lambda(s_\alpha)}I_P^G(s_\alpha\cdot\sigma\otimes\chi^{-1})
\xrightarrow{I_P^G (\varphi_\sigma^{-1}\otimes\id)}I_P^G(\sigma\otimes\chi^{-1})\\
&\xrightarrow{\sim}I_P^G(\overline{\sigma\otimes\chi}^\vee)\xrightarrow{\sim}\overline{I_P^G(\sigma\otimes\chi)}^\vee.
\end{align*}
The last two isomorphisms come from Lemma \ref{lem:inducedrepishermitian}.
Using the map above, we take
\begin{align}\label{eq:hermitianformpoging1}
I_P^G(\sigma\otimes\chi)\times I_P^G(\sigma\otimes\chi)&\xrightarrow{\sim}I_P^G(\sigma\otimes\chi)\times \overline{I_P^G(\sigma\otimes\chi)^\vee}\rightarrow\mathbb{C}\\
\nonumber
(v, v')&\longmapsto \langle v, v'\rangle_{\chi}.
\end{align}
This way, we obtain a sesquilinear form, but it is not necessarily Hermitian. Later, we will do some more work to construct a Hermitian form using (\ref{eq:hermitianformpoging1}).

For now, with (\ref{eq:hermitianformpoging1}), we have a family of nondegenerate sesquilinear forms on the family of representations $I_P^G(\sigma\otimes\chi)$, $\chi\in \Hom(M/Z(G),\mathbb{R}_{>0})$. It only fails when $J'_{\bar{P}|P}(\sigma\otimes\chi)$ is not invertible. The maps depend rationally on $\chi$, as we will prove below.

\begin{lem}\label{lem:hermitianformscontinuous}
The sesquilinear forms constructed in \textnormal{(\ref{eq:hermitianformpoging1})} are $G$-invariant and depend rationally on $\chi\in\Hom(M/Z(G),\mathbb{R}_{>0})$. 
\end{lem}
\begin{proof}

    To show that the sesquilinear forms depend rationally on $\chi$, let us consider each individual map in the composition. It turns out that most of the maps are independent of $\chi$, but let us still describe each of them explicitly.

First of all, we consider the map
\[
J'_{\bar{P}|P}(\sigma\otimes\chi)=(\chi(h_\alpha^\vee)-1)J_{\bar{P}|P}:I_P^G (\sigma\otimes\chi)\rightarrow I_{\bar{P}}^G(\sigma\otimes\chi).
\]
We may identify $I_P^G(\sigma\otimes\chi)\cong I_{K\cap P}^K(\sigma)$ and $I_{\bar{P}}^G\cong I_{K\cap\bar{P}}^K(\sigma)$, to consider it as a map 
\[
J'_{\bar{P}|P}(\sigma\otimes\chi):I_{K\cap P}^K(\sigma)\rightarrow I_{K\cap\bar{P}}^K(\sigma),
\]
and it makes sense to ask if this map is rational in the variable $\chi$, see \cite[\S IV.1]{Wal} for details. By \cite[Th\'eor\`eme IV.1.1]{Wal}, $J_{\bar{P}|P}$ is rational in $\chi$. Moreover $\chi(h_\alpha^\vee)$ is rational in $\chi$, so it follows that indeed $J'_{\bar{P}|P}(\sigma\otimes\chi)=(\chi(h_\alpha^\vee)-1)J_{\bar{P}|P}(\sigma\otimes\chi)$ is rational in $\chi$.

The next step is the map 
\[
\lambda(s_\alpha):I_{\bar{P}\cap K}^K(\sigma)\cong I_{\bar{P}}^G(\sigma\otimes\chi)\rightarrow I_P^G(s_\alpha\cdot \sigma\otimes \chi^{-1})\cong I_{P\cap K}^K (s_\alpha\cdot \sigma).
\]
It is given by $f\mapsto [g\mapsto f(s_\alpha^{-1}g)]$ which is independent of $\chi$, so rational in $\chi$.

Next up is the map 
\[
I_P^G(\varphi_\sigma^{-1}\otimes\id):I_{P\cap K}^K (s_\alpha\cdot\sigma)\cong I_P^G(s_\alpha\cdot \sigma\otimes \chi^{-1})\rightarrow I_P^G(\sigma\otimes\chi^{-1})\cong I_{P\cap K}^K(\sigma).
\]
It is given by $f\mapsto [g\mapsto(\varphi_{\sigma}^{-1}\otimes\id)f(g)]$ which is again independent of $\chi$.

Then we have the map
\[
I_{P\cap K}^K(\sigma)\cong I_P^G(\sigma\otimes\chi^{-1})\xrightarrow{\sim} I_P^G(\overline{\sigma\otimes\chi}^\vee)\cong I_{P\cap K}^K(\overline{\sigma}^\vee).
\]
The isomorphism in the middle comes from unitarity of $\sigma$ and the assumption that $\chi\in\Hom(M,\mathbb{R}_{>0})$, so $\overline{\chi}^\vee\cong \chi^{-1}$. On the level of vector spaces, this map comes from the pairing $V_\sigma\times V_\sigma\rightarrow \mathbb{C}$ inducing the isomorphism $V_\sigma\cong \overline{V_\sigma}^\vee$, which is independent of $\chi$. 

The next map is 
\[
I_{P\cap K}^K (\overline{\sigma}^\vee)\cong I_P^G(\overline{\sigma\otimes\chi}^\vee)\xrightarrow{\sim}I_P^G(\overline{\sigma\otimes\chi})^\vee \cong I_{P\cap K}^K (\overline{\sigma})^\vee,
\]
i.e. the map that describes how taking the contragredient representation commutes with parabolic induction. We have an $M$-invariant pairing $\overline{\sigma\otimes\chi}\times\overline{\sigma\otimes\chi}^\vee\rightarrow\mathbb{C}$. We use it to build a  $K$-invariant pairing
\[
I_{P\cap K}^K(\overline{\sigma\otimes\chi})\times I_{P\cap K}^K(\overline{\sigma\otimes\chi}^\vee)\rightarrow \mathbb{C}, (f_1,f_2)\mapsto \int_K \langle f_1(k),f_2(k)\rangle\td  k,
\]
which gives an isomorphism $I_{P\cap K}^K(\overline{\sigma}^\vee)\cong I_{P\cap K}^K(\overline{\sigma\otimes\chi}^\vee)\cong I_{P\cap K}^K(\overline{\sigma\otimes\chi})^\vee\cong I_{P\cap K}^K (\overline{\sigma})^\vee$ independent of $\chi$. By the Iwasawa decomposition $G = P K$, $\int_K$ is up to 
a positive scalar factor the same as $\int_{P \backslash G}$, hence the pairing is also $G$-invariant, and it induces the isomorphism $I_P^G(\overline{\sigma\otimes\chi}^\vee)\cong I_P^G(\overline{\sigma\otimes\chi})^\vee$, independent of $\chi$.

Next, we need to consider the map
\[
I_{P\cap K}^K(\overline{\sigma})^\vee\cong I_P^G(\overline{\sigma\otimes\chi})^\vee\xrightarrow{\sim}\overline{I_P^G(\sigma\otimes\chi)}^\vee\cong \overline{I_{P\cap K}^K (\sigma)}^\vee,
\]
which shows that taking the complex conjugate representation commutes with parabolic induction. Again, this map is independent of $\chi$. 

Thus, the forms $\langle \ , \ \rangle_\chi$ depend rationally on $\chi$. The sesquilinear forms are $G$-invariant because in the construction in (\ref{eq:hermitianformpoging1}), the first arrow is a composition of homomorphisms of $G$-representations (i.e. the homomorphisms described above), so it is $G$-equivariant, and the second arrow is $G$-invariant by definition. 
\end{proof}

We use the sesquilinear forms in (\ref{eq:hermitianformpoging1}) to build a family of Hermitian forms on the representations $I_P^G(\sigma\otimes\chi)$, which is positive definite on $I_P^G(\sigma)$, and nondegenerate if $\mu(\sigma\otimes\chi)\not=\infty$. We do this by multiplying with a suitable function $f(\chi)$, which depends continuously on $\chi$.

\begin{lem}\label{lem:hermitianformnopole}
    Suppose $\chi\in\Hom(M/Z(G),\mathbb{R}_{>0})$ and $\mu(\sigma\otimes\chi)\not=\infty$. Then $I_P^G(\sigma\otimes\chi)$ is irreducible, and on this representation we have the sesquilinear form as constructed in \textnormal{(\ref{eq:hermitianformpoging1})}. We can multiply it by $f(\chi)\in S^1$, depending continuously on $\chi$, such that
    \begin{equation}\label{eq:hermitianform}
    \langle v, v' \rangle_{I_P^G(\sigma\otimes \chi)}:=f(\chi)\langle v, v' \rangle_{\chi}
    \end{equation}
    is a G-invariant Hermitian form on $I_P^G(\sigma\otimes\chi)$ which is nondegenerate, and positive definite for $\chi=1$.
\end{lem}
\begin{proof}
First we consider the case $\chi=1$. Since $\sigma$ is unitary, $I_P^G(\sigma)$ is also unitary, see \cite[\S IV.2.3]{Ren}. The zero $\mu(\sigma)=0$ corresponds to a pole of order $1$ of $J_{\bar{P}|P}(\sigma\otimes \chi)$ at $\chi=1$. The construction of $J'_{\bar{P}|P}(\sigma\otimes \chi)$ gets rid of the pole at $\chi=1$, therefore $J'_{\bar{P}|P}(\sigma)$ is invertible. Thus, the construction in (\ref{eq:hermitianformpoging1}) gives us a $G$-invariant sesquilinear form $\langle \ , \ \rangle_{\chi=1}$ on $I_P^G(\sigma)$.

By \cite[Proposition 6.3]{tang2017principal}, $I_P^G(\sigma)$ is irreducible. By irreducibility of $I_P^G(\sigma)$, there is up to scalars only one Hermitian form on $I_P^G(\sigma)$, see \cite[112]{Ren}. By unitarity of $I_P^G(\sigma)$, there exists a positive definite Hermitian form $\langle \ , \ \rangle'_{\chi=1}$ on $I_P^G(\sigma)$. Then $\langle \ , \ \rangle_{\chi=1}$ must be a complex scalar multiple of $\langle \ , \ \rangle'_{\chi=1}$, since it is $G$-invariant, sesquilinear and $I_P^G(\sigma)$ is irreducible. We may pick $f(1)\in S^1$ such that $f(1)\langle \ , \ \rangle_{\chi=1}$ equals $\langle \ , \ \rangle'_{\chi=1}$ multiplied by a positive real scalar. In this way, we obtain a Hermitian form on $I_P^G(\sigma)$ which is positive definite.

Now we consider $\chi\not=1$. For $\chi\in \Hom(M/Z(G),\mathbb{R}_{>0}$), we have that $s_\alpha(\chi)=\chi^{-1}$. For $\chi\not=1$, we have $\mu(\sigma\otimes\chi)\not= 0$ by Proposition \ref{prop:zeros} and $s_\alpha(\sigma\otimes\chi)\cong \sigma\otimes\chi^{-1}\not\cong\sigma\otimes\chi$. Hence, by \cite[Proposition 6.2]{tang2017principal}, $I_P^G(\sigma\otimes\chi)$ is irreducible. Since $\mu(\sigma\otimes\chi)\not= 0$ and $\mu(\sigma\otimes\chi)\not= \infty$, we have that $J'_{\bar{P}|P}(\sigma\otimes\chi)$ is invertible. Therefore, we have the $G$-invariant sesquilinear form $\langle \ , \ \rangle_{\chi}$ on $I_P^G(\sigma\otimes\chi)$, as constructed in (\ref{eq:hermitianformpoging1}). 

By Lemma \ref{lem:inducedrepishermitian}, there exists a nondegenerate $G$-invariant Hermitian form $\langle \ , \ \rangle_{\chi}'$ on $I_P^G(\sigma\otimes\chi)$, which is unique up to real scalars. For every $\chi$, fix one such $\langle \ , \ \rangle_{\chi} '$ and scale it, so that we can write 
\[
\langle \ , \ \rangle_\chi =e^{i\phi(\chi)}\langle \ , \ \rangle_{\chi}'.
\]
Then the choice of $\langle \ , \ \rangle_{\chi}'$ is unique up to multiplication by $\pm 1$. Hence, $e^{i\phi(\chi)}$ is also unique up to multiplication by $\pm 1$. Moreover, since $\langle \ , \ \rangle_{\chi}$ depends continuously on $\chi$, 
\[
|\langle \ , \ \rangle_{\chi}|=|e^{i\phi(\chi)}||\langle \ , \ \rangle_{\chi}'|=|\langle \ , \ \rangle_{\chi}'|
\]
depends continuously on $\chi$. Recall that we already fixed $f(1)=e^{-i\phi(1)}$ such that $f(1)\langle \ , \ \rangle_{\chi=1}=\langle \ , \ \rangle_{\chi=1}'$ is positive definite. Consider $0\not =v\in I_{P\cap K}^K(V_\sigma)$, then $\langle v , v \rangle_{\chi=1}'>0$. Consider a neighbourhood of $1$ such that for all $\chi$ in this neighbourhood, $\langle v,v \rangle_{\chi}'\not =0$. This is possible since $|\langle \ , \ \rangle_{\chi}'|$ depends continuously on $\chi$. The $\langle \ , \ \rangle_{\chi}'$ were chosen uniquely up to a sign, so it is possible to choose them in such a way that $\langle v , v \rangle_{\chi}'>0$ for all $\chi$ in this neighbourhood of $1$. In this way, $\langle \ , \ \rangle_{\chi}'$ depends continuously on $\chi$ in this neighbourhood. In fact, for all $\chi$ we can pick the sign of $\langle \ , \ \rangle_{\chi}'$ such that the Hermitian forms will depend continuously on $\chi$. Then also $f(\chi)=e^{-i\phi(\chi)}$ depends continuously on $\chi$, and we are done.
\end{proof}

\begin{lem}\label{lem:hermitianformrestricts}
    Let $(\pi,V)$ be a smooth Hermitian representation of $G$. Let $C$ be a compact open subgroup of $G$ and let $V^C$ be the space of $C$-fixed vectors in $V$. Then the nondegenerate $G$-invariant Hermitian form on $V$ restricts to a nondegenerate $C$-invariant Hermitian form on $V^C$. 
\end{lem}
\begin{proof}
Since $C$ is compact, we can write 
\[
V=V^C\oplus V_C
\]
where $V_C$ is spanned by the vectors $c\cdot v-v$ for $c\in C$, $v\in V$. Since $(\pi, V)$ is Hermitian, we have a $G$-invariant nondegenerate Hermitian form on $V$, \[\langle \ \cdot \ , \  \cdot \ \rangle:V\times V=(V^C\oplus V_C)\times (V^C\oplus V_C)\rightarrow \mathbb{C}.\]
Consider the projection operator 
\[
p_C:=\frac{1}{|C|}\int_C \pi(c) \td c=\begin{cases}
    1 \textnormal{ on } V^C\\
    0 \textnormal{ on } V_C.
\end{cases}
\]
Since the Hermitian form on $V$ is $G$-invariant, in particular $C$-invariant, we have that 
\[
\langle c\cdot v,v'\rangle=\langle v, c^{-1}\cdot v'\rangle
\]
for all $v,v'\in V$, $c\in C$. It follows that for all $v\in V^C$, $v'\in V_C$, 
\[
\langle v, v'\rangle =\langle p_C v, v'\rangle =\langle v, p_C v'\rangle =\langle v, 0 \rangle =0.
\] So the pairing between $V^C$ and $V_C$ is zero, but the form on the whole space $V$ is nondegenerate, therefore the form restricted to $V^C$ must be nondegenerate.
\end{proof}

For a Hermitian form on a finite dimensional vector space $V$, we can consider its signature $(\dim(V^+) ,\dim(V^0), \dim(V^-))$, where $V=V^+ \oplus V^0\oplus V^-$ and the form is positive definite on $V^+$, zero on $V^0$, negative definite on $V^-$. For Hermitian forms on admissible representations, we can consider the signature by restricting to open compact subgroups, which will be done to prove the next result.

\begin{lem}\label{lem:unitary}
    Suppose $\mu(\sigma)=0$, $\chi\in\Hom(M/Z(G),\mathbb{R}_{>0})$. If $\mu(\sigma\otimes\chi^t)$ does not have a pole for $t\in[0,\tau]$, $\tau\in\mathbb{R}_{>0}$, then $I_P^G(\sigma\otimes\chi^\tau)$ is unitary. 
\end{lem}
\begin{proof}
    By Lemma \ref{lem:hermitianformnopole}, for each $t\in[0,\tau]$, the representation $I_P^G(\sigma\otimes\chi^t)$ is Hermitian, and we have the Hermitian form as constructed in (\ref{eq:hermitianform}). The representations are realized on the same vector space, say $I_{P\cap K}^K (V_\sigma)$. By Lemma \ref{lem:hermitianformnopole}, the Hermitian form on $I_P^G(\sigma)$ is positive definite, i.e. $\langle f, f\rangle_{I_P^G(\sigma)}>0$. The function $t\mapsto\langle f, f\rangle_{I_P^G(\sigma\otimes\chi^t)}$ is continuous for $f\in I_{P\cap K}^K (V_\sigma)$, see Lemma \ref{lem:hermitianformscontinuous} and Lemma \ref{lem:hermitianformnopole}. Let $C\subset K$ be a compact open subgroup fixing $f$. Recall that $I_P^G(\sigma\otimes\chi^t)$ is irreducible, hence admissible by \cite[Theorem 6.3]{FratilaPrasad}, therefore $\dim I_{P\cap K}^K (V_\sigma)^C$ is finite. By Lemma \ref{lem:hermitianformrestricts}, the $G$-invariant nondegenerate Hermitian form on $I_P^G(\sigma\otimes\chi^t)$ constructed in (\ref{eq:hermitianform}) restricts to a $C$-invariant nondegenerate Hermitian form on $I_{P\cap K}^K(V_\sigma)^C$. We have a continuous family of such forms parametrized by $t\in[0,\tau]$. For $t=0$, recall that we have a positive definite form, so the form has signature $(\dim I_{P\cap K}^K(V_\sigma)^C,0,0)$. For any $t\in [0,\tau]$, the form is nondegenerate, so with signature $(p_t,0,q_t)$ with $p_t+q_t=\dim I_{P\cap K}^K(V_\sigma)^C$. If for some $t$ we have $q_t>0$, then since the forms vary continuously, we must have that for some $t'$ the form is degenerate, a contradiction. So all the Hermitian forms on $I_{P\cap K}^K(V_\sigma)^C$ for $t\in[0,\tau]$ are positive definite. We can do this for all $f$ and $C$ fixing $f$, hence all the representations $I_P^G(\sigma\otimes\chi^t)$ for $t\in [0,\tau]$ are unitary.
\end{proof}

We have just shown that if $\mu(\sigma\otimes \chi^t)$ does not have a pole for any $t\in\mathbb{R}_{\geq0}$, then all representations $I_P^G(\sigma\otimes\chi^t)$ are unitary. We want to show that this is not possible, so that we obtain a pole for some $t$.

For the next result, we consider the maximal $C^*$-algebra $C^*(G)$ of $G$, which is a completion of $\mathcal{H}(G)$ such that the irreducible smooth unitary representations of $G$ are in bijection with the irreducible representations of $C^*(G)$. This bijection is given by
\[
\Irr(G)^{\unit}\rightarrow \Irr(C^*(G)), (\pi, V)\mapsto (\pi, \mathcal{V}),
\]
where $\mathcal{V}$ denotes the Hilbert space completion of $V$, and the action of $C^*(G)$ on $\mathcal{V}$ is (for functions $f\in\mathcal{H}(G)$) given by
\begin{equation}\label{eqn:repC*algebra}
\pi(f)(v)=\int_Gf(g)\pi(g)v\td g.
\end{equation}
Every unitary representation $(\pi,V)$ of $G$ gives a homomorphism of $C^*$-algebras $\pi:C^*(G)\rightarrow \mathcal{B}(\mathcal{V})$, where $\mathcal{B}(\mathcal{V})$ is the $C^*$-algebra of bounded linear operators on $\mathcal{V}$. Since it is a homomorphism of $C^*$-algebras, it must be norm-decreasing, see \cite[Theorem 2.1.7]{Murphy}.

\begin{lem}\label{lem:mupole}
    Suppose $\mu(\sigma)=0$ and let $\chi\in\Hom(M/Z(G),\mathbb{R}_{>0})$ such that $\chi(h_\alpha^\vee)>1$. Then $\mu(\sigma\otimes \chi^t)$ has a pole for some $t\in\mathbb{R}_{>0}$.
\end{lem}
\begin{proof}
    Suppose not, then by Lemma \ref{lem:unitary}, $\pi^t:= I_P^G(\sigma\otimes\chi^t)$ is a unitary $G$-representation for all $t\in \mathbb{R}_{\geq 0}$. 
    Let $v\in I_{P\cap K}^K (V_\sigma)\cong I_P^G(V_\sigma)$ such that $\supp(v)= (P\cap K)\cdot U_{-\alpha,r}$ for some $r$, and such that $v$ is right $U_{-\alpha,r}$-invariant. Then $v(\bar{u})=v(1)$ for all $\bar{u}\in U_{-\alpha,r}$. Let $C$ be a compact open subgroup of $G$ fixing $v$ (this is possible because $I_P^G(V_\sigma)$ is smooth) and put $f=1_{h_\alpha^\vee C}$. Then $f\in\mathcal{H}(G)\subset C^*(G)$, where $C^*(G)$ is the maximal $C^*$-algebra of $G$, so $||f||_{C^*(G)}\in\mathbb{R}_{\geq0}$. We can view $\pi^t$ as a representation of $C^*(G)$, as written in (\ref{eqn:repC*algebra}). We have
    \begin{align*}
      \pi^t(f)v&= \int_Gf(g)\pi^t(g)v\td g=\int_{h_\alpha^\vee C}\pi^t(g)v\td g= \int_{h_\alpha^\vee C}(\sigma\otimes\chi^t)(g)v\td g\\
      &=\int_{h_\alpha^\vee C}(\sigma\otimes \chi^t)(h_\alpha^\vee)v \td g=\mr{vol}(h_\alpha^\vee C)\cdot(\sigma\otimes\chi^t)(h_\alpha^\vee)v,
    \end{align*}
    where the second to last equality holds because $C$ is a compact open subgroup fixing $v$, so the only nontrivial contribution of the action of $h_\alpha^\vee C$ on $v$ comes from $h_\alpha^\vee$. Note that $\pi^t(f)v$ is again right $U_{-\alpha,r}$-invariant, and supported in $(P\cap K)\cdot U_{-\alpha,r}$.
    
    Since $\pi^t(f)v$ lies in $I_{P\cap K}^K(V_\sigma)$, it also lies in the Hilbert space completion $\mathcal{V}$ of $I_{P\cap K}^K (V_\sigma)$. Let $\mathcal{B}(\mathcal{V})$ be the $C^*$-algebra of bounded linear operators on $\mathcal{V}$, with the operator norm denoted by $||\, . \, ||_\op$. 
    The operator norm of $\pi^t(f)$ is given by
    \[
    ||\pi^t(f)||_{\op}=\sup _{v'\not=0}\frac{||\pi^t(f)v'||_{\mathcal{V}}}{||v'||_{\mathcal{V}}}, \textnormal{ where } ||v'||_{\mathcal{V}}^2=\int_K|| v'(k)||^2\td k.
    \]
    Let $v$ with $\supp(v)=(P\cap K)\cdot U_{-\alpha,r}$ as above. Then 
    \begin{align*}
    ||v||_{\mathcal{V}}^2&=\int_K||v(k)||^2\td k = \int_{(P\cap K)\cdot U_{-\alpha,r}}||v(k)||^2 \td k =\mr{vol}(U_{-\alpha,r})\int_{P\cap K}|| v(k)||^2\td k\\
    &=\mr{vol}(U_{-\alpha,r})\int_{P\cap K}|| \sigma(k)v(1)||^2 \td k=\mr{vol}(U_{-\alpha,r})\mr{vol}(P\cap K)|| v(1)||^2,
    \end{align*}
    since $v\in I_{P\cap K}^K(V_\sigma)$ is right $U_{-\alpha,r}$-invariant and $\sigma$ is unitary.
    Similarly, we have 
    \begin{align*}
    ||\pi^t(f)v||_{\mathcal{V}}^2&=\int_K|| \pi^t(f)v(k)||^2\td k=\int_{(P\cap K)\cdot U_{-\alpha,r}} || \pi^t(f)v(k)||^2\td k\\
    &=\mr{vol}(U_{-\alpha,r})\mr{vol}(h_\alpha^\vee C)^2\int_{P\cap K}||(\sigma\otimes\chi^t)(h_\alpha^\vee) v(k)||^2 \td k\\
    &=\mr{vol}(U_{-\alpha,r})\mr{vol}(h_\alpha^\vee C)^2(\chi^t(h_\alpha^\vee))^2\int_{P\cap K}||\sigma(h_\alpha^\vee)\sigma(k)v(1)||^2\td k\\
    &=\mr{vol}(U_{-\alpha,r})\mr{vol}(h_\alpha^\vee C)^2(\chi^t(h_\alpha^\vee))^2\mr{vol}(P\cap K)||v(1)||^2\\
    &=\mr{vol}(h_\alpha^\vee C)^2(\chi^t(h_\alpha^\vee))^2||v||_\mathcal{V}^2.
    \end{align*} 
    Recall that since $\pi^t$ is unitary, it gives a homomorphism of $C^*$-algebras $\pi^t:C^*(G)\rightarrow \mathcal{B}(\mathcal{V})$, which is norm-decreasing.  Since $\chi(h_\alpha^\vee)>1$, we have for $t\rightarrow\infty$ that
    \[
    ||\pi^t(f)v||_\mathcal{V}/||v||_\mathcal{V}=\mr{vol}(h_\alpha^\vee C)\chi^t(h_\alpha^\vee)\rightarrow \infty \textnormal{, hence } \lim_{t\rightarrow\infty}||\pi^t(f)||_{\op}=\infty
    \]
    which contradicts $||\pi^t(f)||_{\op}\leq ||f||_{C^*(G)}$. So $\mu(\sigma\otimes\chi^t)=\infty$ for some $t>0$.
\end{proof}

Now we are ready to prove the main result. We drop the assumption $\sigma=\sigma_0$ that was used to prove the results in this subsection so far.

\begin{thm} There exists a unitary $\sigma_0\in \Irr(M)_{[M,\sigma]}$ such that 
\[
\mu(\sigma_0\otimes \chi)=\mu(\sigma_0,z) =c_{\mu}\cdot \frac{(1-z)(1-z^{-1})}{(1-qz)(1-qz^{-1})}\frac{(1+z)(1+z^{-1})}{(1+q'z)(1+q'z^{-1})},
\]
where $c_{\mu}\in\mathbb{R}_{>0}$ and $q\geq1, \ q'\geq1$.
\end{thm}
\begin{proof}
Suppose that $\mu(\sigma,z)$ is constant, then we can take $q=q'=1$. Suppose that $\mu(\sigma,z)$ is not constant, then there exists a unitary $\sigma_0\in \Irr(M)_{[M,\sigma]}$ as in Proposition \ref{prop:zeros}. By Lemma \ref{lem:formmu}, $\mu(\sigma_0,z)$ has the form
    \[
    \mu(\sigma_0,z)=c_{\mu}\prod_i (z-\lambda_i)^{n_i}(z^{-1}-\overline{\lambda_i})^{n_i},
    \]
 where $\lambda_i\in\mathbb{C}$, $n_i\in \mathbb{Z}$, $c_{\mu}\in\mathbb{R}_{>0}$. By Proposition \ref{prop:muatinfinity}, $\mu(\sigma_0,z)$ is not zero at $z=0$ and $z=\infty$. By Proposition \ref{prop:zeros}, the zero at $z=1$ is of order $2$. So $\mu(\sigma_0,z)$ has a factor $(1-z)(1-z^{-1})$ in the numerator. By Lemma \ref{lem:mupole}, $\mu(\sigma_0)=0$ implies that $\mu(\sigma_0\otimes\chi^t)$ has a pole for some $t\in \mathbb{R}_{>0}$, $\chi\in \Hom(M/Z(G),\mathbb{R}_{>0})$ with $\chi(h_\alpha^\vee)>1$. So $\mu(\sigma_0,z)$ has a pole at $q=\chi^t(h_\alpha^\vee)\in\mathbb{R}_{>1}$ and $\mu$ has a factor $(1-qz)(1-qz^{-1})$ in the denominator. The only other possible zero of $\mu(\sigma_0,z)$ is at $z=-1$, and again of multiplicity $2$, so $\mu$ has a factor $(1+z)(1+z^{-1})$ in the numerator (if $\mu$ is not zero at $z=-1$, we can cancel out this factor by taking $q'=1$ in the formula). Suppose that $\mu(\sigma_0,-1)=0$ and let $\chi'$ be a unitary character such that $z=\chi'(h_\alpha^\vee)=-1$. Then $\mu(\sigma_0\otimes \chi')=0$, and by Lemma \ref{lem:mupole}, $\mu((\sigma_0\otimes\chi')\otimes\chi^t)=\mu(\sigma_0\otimes\chi'\chi^t)$ has a pole for some $\chi\in\Hom(M/Z(G),\mathbb{R}_{>0})$, $t\in\mathbb{R}_{>0}$. This implies that $\mu(\sigma_0,z)$ has a pole at $z=\chi'(h_\alpha^\vee)\chi^t(h_\alpha^\vee)=-1\cdot q'$ with $q'\in\mathbb{R}_{>1}$. So $\mu$ has a factor $(1+q'z)(1+q'z^{-1})$ in the denominator.

 It remains to show that $\mu(\sigma_0,z)$ does not have any additional factors. More factors in the numerator would mean that $\mu(\sigma_0,z)$ has extra zeros in $\mathbb{C}^\times$, which is not possible. More factors in the denominator would mean that $\mu(\sigma_0,\infty)=\mu(\sigma_0,0)=0$, which is not possible by Proposition \ref{prop:muatinfinity}.
\end{proof}

\begingroup
\setlength{\emergencystretch}{8em}
\printbibliography

@article {Wal,
    AUTHOR = {Waldspurger, J.-L.},
     TITLE = {La formule de {P}lancherel pour les groupes {$p$}-adiques
              (d'apr\`es {H}arish-{C}handra)},
   JOURNAL = {J. Inst. Math. Jussieu},
  FJOURNAL = {Journal of the Institute of Mathematics of Jussieu. JIMJ.
              Journal de l'Institut de Math\'{e}matiques de Jussieu},
    VOLUME = {2},
      YEAR = {2003},
    NUMBER = {2},
     PAGES = {235--333},
      ISSN = {1474-7480,1475-3030},
   MRCLASS = {22E35 (22E50)},
  MRNUMBER = {1989693},
MRREVIEWER = {Rebecca\ Herb},
       DOI = {10.1017/S1474748003000082},
       URL = {https://doi.org/10.1017/S1474748003000082},
}

@article {Sil,
    AUTHOR = {Silberger, Allan J.},
     TITLE = {Special representations of reductive {$p$}-adic groups are not
              integrable},
   JOURNAL = {Ann. of Math. (2)},
  FJOURNAL = {Annals of Mathematics. Second Series},
    VOLUME = {111},
      YEAR = {1980},
    NUMBER = {3},
     PAGES = {571--587},
      ISSN = {0003-486X},
   MRCLASS = {22E50 (22E35)},
  MRNUMBER = {577138},
MRREVIEWER = {Joe\ Repka},
       DOI = {10.2307/1971110},
       URL = {https://doi.org/10.2307/1971110},
}

@book {Ren,
    AUTHOR = {Renard, David},
     TITLE = {Repr\'{e}sentations des groupes r\'{e}ductifs {$p$}-adiques},
    SERIES = {Cours Sp\'{e}cialis\'{e}s [Specialized Courses]},
    VOLUME = {17},
 PUBLISHER = {Soci\'{e}t\'{e} Math\'{e}matique de France, Paris},
      YEAR = {2010},
     PAGES = {vi+332},
      ISBN = {978-2-85629-278-5},
   MRCLASS = {22E50 (20G25)},
  MRNUMBER = {2567785},
MRREVIEWER = {Anne-Marie\ H.\ Aubert},
}

@misc{tang2017principal,
      title={Principal series representations of metaplectic groups}, 
      author={Shiang Tang},
      year={2017},
      eprint={1706.05145},
      archivePrefix={arXiv},
      primaryClass={math.RT}
}

@misc{solleveld2023endomorphismalgebrasheckealgebras,
      title={Endomorphism algebras and Hecke algebras for reductive p-adic groups}, 
      author={Maarten Solleveld},
      year={2023},
      eprint={2005.07899},
      archivePrefix={arXiv},
      primaryClass={math.RT},
      url={https://arxiv.org/abs/2005.07899}, 
}

@article {FratilaPrasad,
    AUTHOR = {Fratila, Dragos and Prasad, Dipendra},
     TITLE = {Homological duality for covering groups of reductive
              {$p$}-adic groups},
   JOURNAL = {Pure Appl. Math. Q.},
  FJOURNAL = {Pure and Applied Mathematics Quarterly},
    VOLUME = {18},
      YEAR = {2022},
    NUMBER = {5},
     PAGES = {1867--1950},
      ISSN = {1558-8599,1558-8602},
   MRCLASS = {11F70 (22E55)},
  MRNUMBER = {4538041},
MRREVIEWER = {Kimball\ L.\ Martin},
}

@book {CGP,
    AUTHOR = {Conrad, Brian and Gabber, Ofer and Prasad, Gopal},
     TITLE = {Pseudo-reductive groups},
    SERIES = {New Mathematical Monographs},
    VOLUME = {26},
   EDITION = {Second},
 PUBLISHER = {Cambridge University Press, Cambridge},
      YEAR = {2015},
     PAGES = {xxiv+665},
      ISBN = {978-1-107-08723-1},
   MRCLASS = {20G15 (14L15)},
  MRNUMBER = {3362817},
       DOI = {10.1017/CBO9781316092439},
       URL = {https://doi.org/10.1017/CBO9781316092439},
}

@book {KaPr,
    AUTHOR = {Kaletha, Tasho and Prasad, Gopal},
     TITLE = {Bruhat-{T}its theory---a new approach},
    SERIES = {New Mathematical Monographs},
    VOLUME = {44},
 PUBLISHER = {Cambridge University Press, Cambridge},
      YEAR = {2023},
     PAGES = {xxx+718},
      ISBN = {978-1-108-83196-3},
   MRCLASS = {20E42 (11F70 20G25 22E50)},
  MRNUMBER = {4520154},
MRREVIEWER = {Corina\ Ciobotaru},
}

@article {Sau,
    AUTHOR = {Sauvageot, Francois},
     TITLE = {Principe de densit\'e{} pour les groupes r\'eductifs},
   JOURNAL = {Compositio Math.},
  FJOURNAL = {Compositio Mathematica},
    VOLUME = {108},
      YEAR = {1997},
    NUMBER = {2},
     PAGES = {151--184},
      ISSN = {0010-437X,1570-5846},
   MRCLASS = {22E50 (22E55 43A05 43A70)},
  MRNUMBER = {1468833},
MRREVIEWER = {Volker\ J.\ Heiermann},
       DOI = {10.1023/A:1000216412619},
       URL = {https://doi.org/10.1023/A:1000216412619},
}

@book {Murphy,
    AUTHOR = {Murphy, Gerard J.},
     TITLE = {{$C^*$}-algebras and operator theory},
 PUBLISHER = {Academic Press Inc., Boston, MA},
      YEAR = {1990},
     PAGES = {x+286},
      ISBN = {0-12-511360-9},
   MRCLASS = {46Lxx (46-01)},
  MRNUMBER = {1074574},
MRREVIEWER = {E.\ Gerlach},
}

@article {WenWeiLi,
    AUTHOR = {Li, Wen-Wei},
     TITLE = {La formule des traces pour les rev\^etements de groupes
              r\'eductifs connexes. {II}. {A}nalyse harmonique locale},
   JOURNAL = {Ann. Sci. \'Ec. Norm. Sup\'er. (4)},
  FJOURNAL = {Annales Scientifiques de l'\'Ecole Normale Sup\'erieure.
              Quatri\`eme S\'erie},
    VOLUME = {45},
      YEAR = {2012},
    NUMBER = {5},
     PAGES = {787--859},
      ISSN = {0012-9593,1873-2151},
   MRCLASS = {43A60 (11F70 11F72)},
  MRNUMBER = {3053009},
MRREVIEWER = {Joseph\ Max\ Rosenblatt},
       DOI = {10.24033/asens.2178},
       URL = {https://doi.org/10.24033/asens.2178},
}

@article {endomorphismalgebras,
    AUTHOR = {Solleveld, Maarten},
     TITLE = {Endomorphism algebras and {H}ecke algebras for reductive
              {$p$}-adic groups},
   JOURNAL = {J. Algebra},
  FJOURNAL = {Journal of Algebra},
    VOLUME = {606},
      YEAR = {2022},
     PAGES = {371--470},
      ISSN = {0021-8693,1090-266X},
   MRCLASS = {22E50 (20C08 20G25)},
  MRNUMBER = {4432237},
MRREVIEWER = {Daniel\ Goldstein},
       DOI = {10.1016/j.jalgebra.2022.05.017},
       URL = {https://doi.org/10.1016/j.jalgebra.2022.05.017},
}

@article {silbergerknappstein,
    AUTHOR = {Silberger, Allan J.},
     TITLE = {The {K}napp-{S}tein dimension theorem for {$p$}-adic groups},
   JOURNAL = {Proc. Amer. Math. Soc.},
  FJOURNAL = {Proceedings of the American Mathematical Society},
    VOLUME = {68},
      YEAR = {1978},
    NUMBER = {2},
     PAGES = {243--246},
      ISSN = {0002-9939,1088-6826},
   MRCLASS = {22E50},
  MRNUMBER = {492091},
       DOI = {10.2307/2041781},
       URL = {https://doi.org/10.2307/2041781},
}

@article {BrTi,
    AUTHOR = {Bruhat, F. and Tits, J.},
     TITLE = {Groupes r\'eductifs sur un corps local: I. Donn\'ees radicielles valu\'ees},
   JOURNAL = {Inst. Hautes \'Etudes Sci. Publ. Math.},
  FJOURNAL = {Institut des Hautes \'Etudes Scientifiques. Publications
              Math\'ematiques},
    NUMBER = {41},
      YEAR = {1972},
     PAGES = {5--251},
      ISSN = {0073-8301,1618-1913},
   MRCLASS = {20G25 (22E20)},
  MRNUMBER = {327923},
MRREVIEWER = {M.\ Krusemeyer},
       URL = {http://www.numdam.org/item?id=PMIHES_1972__41__5_0},
}

@article {BerZelInducedRep,
	AUTHOR = {Bernstein, I. N. and Zelevinsky, A. V.},
	TITLE = {Induced representations of reductive {$p$}-adic groups. {I}},
	JOURNAL = {Ann. Sci. \'Ecole Norm. Sup. (4)},
	FJOURNAL = {Annales Scientifiques de l'\'Ecole Normale Sup\'erieure.
	Quatri\`eme S\'erie},
	VOLUME = {10},
	YEAR = {1977},
	NUMBER = {4},
	PAGES = {441--472},
	ISSN = {0012-9593},
	MRCLASS = {22E50},
	MRNUMBER = {579172},
	URL = {http://www.numdam.org/item?id=ASENS_1977_4_10_4_441_0},
}

@article {parameters,
	AUTHOR = {Solleveld, Maarten},
	TITLE = {Parameters of {H}ecke algebras for {B}ernstein components of
	{$p$}-adic groups},
	JOURNAL = {Indag. Math. (N.S.)},
	FJOURNAL = {Koninklijke Nederlandse Akademie van Wetenschappen.
	Indagationes Mathematicae. New Series},
	VOLUME = {36},
	YEAR = {2025},
	NUMBER = {1},
	PAGES = {124--170},
	ISSN = {0019-3577,1872-6100},
	MRCLASS = {22E50 (20C08 20G25)},
	MRNUMBER = {4847675},
	DOI = {10.1016/j.indag.2024.04.005},
	URL = {https://doi.org/10.1016/j.indag.2024.04.005},
}

@article {isogenyrestrictions,
	AUTHOR = {Silberger, Allan J.},
	TITLE = {Isogeny restrictions of irreducible admissible representations
	are finite direct sums of irreducible admissible
	representations},
	JOURNAL = {Proc. Amer. Math. Soc.},
	FJOURNAL = {Proceedings of the American Mathematical Society},
	VOLUME = {73},
	YEAR = {1979},
	NUMBER = {2},
	PAGES = {263--264},
	ISSN = {0002-9939,1088-6826},
	MRCLASS = {22E50},
	MRNUMBER = {516475},
	MRREVIEWER = {G.\ I.\ Ol\cprime shanski\u i},
	DOI = {10.2307/2042302},
	URL = {https://doi.org/10.2307/2042302},
}
\endgroup

\end{document}